\documentclass{amsart}
\usepackage[latin1]{inputenc}
\usepackage{amssymb, enumerate, xspace, graphicx,stmaryrd}
\usepackage[all]{xy}
\SelectTips{cm}{}

\numberwithin{equation}{section}

\numberwithin{subsection}{section}

\allowdisplaybreaks[1]

\newenvironment{enumeratea}
{\begin{enumerate}[\upshape (a)]}
{\end{enumerate}}

\newtheorem*{namedtheorem}{\theoremname}
\newcommand{\theoremname}{testing}

\newtheorem{theorem}{Theorem}[section]
\newtheorem{proposition}[theorem]{Proposition}
\newtheorem{proposition-definition}[theorem]
{Proposition-Definition}
\newtheorem{corollary}[theorem]{Corollary}
\newtheorem{lemma}[theorem]{Lemma}

\theoremstyle{definition}
\newtheorem{definition}[theorem]{Definition}

\newtheorem{remark}[theorem]{Remark}

\theoremstyle{remark}

\newcommand\nome{testing}
\newcommand\call[1]{\label{#1}\renewcommand\nome{#1}}
\newcommand\itemref[1]{\item\label{\nome;#1}}
\newcommand\refall[2]{\ref{#1}~(\ref{#1;#2})}
\newcommand\refpart[2]{(\ref{#1;#2})}

\newcommand\cA{\mathcal{A}} \newcommand\cB{\mathcal{B}}
 
\newcommand\cE{\mathcal{E}} \newcommand\cF{\mathcal{F}}
 \newcommand\cH{\mathcal{H}}
 
 \newcommand\cL{\mathcal{L}}
\newcommand\cM{\mathcal{M}} 
\newcommand\cO{\mathcal{O}} \newcommand\cP{\mathcal{P}}
\newcommand\cQ{\mathcal{Q}} 
\newcommand\cS{\mathcal{S}} \newcommand\cT{\mathcal{T}}
 \newcommand\cV{\mathcal{V}}
\newcommand\cW{\mathcal{W}}

\renewcommand\AA{\mathbb{A}} 
\newcommand\CC{\mathbb{C}} 
 
\newcommand\GG{\mathbb{G}} \newcommand\HH{\mathbb{H}}

 \newcommand\PP{\mathbb{P}}

 \newcommand\ZZ{\mathbb{Z}}

 \newcommand\bT{\mathbf{T}}

\newcommand\rA{\mathrm{A}}

\newcommand\rI{\mathrm{I}} \newcommand\rJ{\mathrm{J}}
 
\newcommand\rM{\mathrm{M}}

 \newcommand\rT{\mathrm{T}}

\newcommand\rmc{\mathrm{c}}

 \newcommand\frg{\mathfrak{g}}

\newcommand\arr{\ifinner\to\else\longrightarrow\fi}

\newcommand\arrto{\ifinner\mapsto\else\longmapsto\fi}
\newcommand\larr{\longrightarrow}
\newcommand\larrto{\longmapsto}

\newcommand{\xarr}{\xrightarrow}

\renewcommand\H{\operatorname{H}}

\newcommand\eqdef{\overset{\mathrm{\scriptscriptstyle def}} =}

\newcommand\into{\hookrightarrow}

\renewcommand\th{^\text{th}}

\def\displaytimes_#1{\mathrel{\mathop{\times}\limits_{#1}}}

\def\displayotimes_#1{\mathrel{\mathop{\bigotimes}\limits_{#1}}}

\renewcommand\hom{\operatorname{Hom}}

\newcommand\pic{\operatorname{Pic}}

\newcommand\spec{\operatorname{Spec}}

\newcommand\generate[1]{\langle #1 \rangle}

\newcommand\id{\mathrm{id}}

\newcommand\coker{\operatorname{coker}}

\newcommand\ddash{\nobreakdash--\hspace{0pt}}

\newdir{ >}{{}*!/-5pt/@{>}}

\newcommand{\catsch}[1]{(\mathrm{Sch}/#1)}

\newcommand\doublelong[2]{\mathbin{\xymatrix{{}\ar@<3pt>[r]^{#1}
\ar@<-3pt>[r]_{#2}&}}}

\newcommand{\underhom}
{\mathop{\underline{\mathrm{Hom}}}\nolimits}

\newcommand{\underaut}
{\mathop{\underline{\mathrm{Aut}}}\nolimits}

\newlength{\ignora}
\newcommand{\hsmash}[1]{\settowidth{\ignora}{#1}#1\hspace{-\ignora}}

\renewcommand{\setminus}{\smallsetminus}

\newcommand{\tr}{\operatorname{tr}}

\newcommand{\T}[1][g]{\cT_{#1}}
\newcommand{\Tbig}[1][g]{\widetilde{\cT}_{#1}}

\newcommand{\OT}[1][g]{\widehat{\cT}_{#1}}

\newcommand{\tsch}{Tschirnhausen\xspace}

\newcommand{\vectp}[1][r,d]{\overline{\cV}_{#1}}
\newcommand{\vect}[1][r,d]{\cV_{#1}}

\newcommand{\vectalt}[1][r,d]{\cV'_{#1}}

\newcommand{\vectaux}[1][r,d]{\cW_{#1}}

\newcommand{\trip}{\operatorname{Trip}}

\newcommand{\pone}{\PP^{1}}

\newcommand{\om}[1][r,d]{\Omega_{#1}}
\newcommand{\mm}[1][r,d]{\rM_{#1}}

\newcommand{\grp}[1][r,d]{\Gamma_{#1}}

\newcommand{\GL}{\mathrm{GL}}

\newcommand{\PGL}{\mathrm{PGL}}

\newcommand{\gm}{\mathbb{G}_{\mathrm{m}}}

\newcommand{\gms}[1]{\mathbb{G}_{\mathrm{m}, #1}}

\newcommand{\sym}{\operatorname{Sym}}
\newcommand{\symthree}[2][]{\sym^{3}_{{#1}}{#2} \otimes_{{#1}} \det\nolimits_{{#1}} {#2}\dual}

\newcommand{\GLe}{\GL_{2, \epsilon}}

\newcommand{\trans}{^{\mathrm{t}}}

\newcommand{\dual}{^{\vee}}
\newcommand\Hg{H_{\GL_2}}
\newcommand\Ve{V_{\epsilon}}
\newcommand{\ch}{\operatorname{c}}

\begin{document}

\title{Stacks of trigonal curves}

\author[Bolognesi]{Michele Bolognesi}
\author[Vistoli]{Angelo Vistoli}

\address[Bolognesi]{Dipartimento di Matematica\\Università di Roma Tre\\Largo San Leonardo Murialdo, 1
00146 Roma\\Italy}
\email{michele.bolognesi@sns.it}

\address[Vistoli]{Scuola Normale Superiore\\Piazza dei Cavalieri 7\\
56126 Pisa\\Italy}
\email{angelo.vistoli@sns.it}



\begin{abstract}
In this paper we study the stack $\T$ of smooth triple covers of a conic; when $g \geq 5$ this stack is embedded $\cM_{g}$ as the locus of trigonal curves. We show that $\T$ is a quotient $[U_{g}/\Gamma_{g}]$, where $\Gamma_g$ is a certain algebraic group and $U_g$ is an open subscheme of a $\Gamma_g$-equivariant vector bundle over an open subscheme of a representation of $\Gamma_g$. Using this, we compute the integral Picard group of $\cT_g$ when $g > 1$. The main tools are a result of Miranda that describes a flat finite triple cover of a scheme $S$ as given by a locally free sheaf $E$ of rank two on $S$, with a section of $\symthree{E}$, and a new description of the stack of globally generated locally free sheaves of  fixed rank and degree on a projective line as a quotient stack. 
\end{abstract}

\maketitle

\section{Introduction}

In moduli theory, stacks have often been constructed as quotients by group actions. If an algebraic stack $\cF$ is a quotient stack $[X/G]$, where $G$ is an algebraic group acting on an algebraic variety $X$, the geometry of $\cF$ is the geometry of the action of $G$ on $X$, and one can apply to $\cF$ the powerful techniques that have been developed for studying invariants of group actions in algebraic topology and algebraic geometry.

Even just knowing that such a presentation exists, even without an explicit description of the action, can be useful: but it is even better when the variety $X$ and the group $G$ are fairly simple, so that this description may be used directly to study $\cF$. This does not seem possible in many cases: for example, the stack $\cM_{g}$ of smooth curves of genus~$g$ is of the form $[X/G]$, but when $g$ is large the space $X$ is complicated, and to our knowledge no general result about $\cM_{g}$ has been proved by exploiting this presentation. (Of course, in Teichmüller theory one represents $\cM_{g}$ as a quotient of an action of the Teichmüller group, which is an infinite discrete group, acting on a ball in $\CC^{3g-3}$, but this description is topological, and it is hard to use it directly to prove algebraic geometric results about $\cM_{g}$.)

It has long been known that in characteristic different from $2$ and $3$, the stack $\cM_{1,1}$ of elliptic curves is a quotient $[(X/\gm)]$, where $X$ is the complement of the hypersurface $4x^{3} + 27y^{2} = 0$ in $\AA^{2}$, and $\gm$ acts with weights $4$ and $6$. This gives an easy proof of the fact, due to Mumford (\cite{mumford-picard}), that the Picard group of $\cM_{1,1}$ is cyclic of order 12. In \cite{vistoli98}, the second author gives a presentation of $\cM_{2}$ as a quotient $[X/\GL_{2}]$, where $X$ is the scheme of smooth binary forms of degree~$6$ in two variables (the action of $\GL_{2}$ on $X$ is a twist of the customary one). As an application he computes the Chow ring of $\cM_{2}$. This was generalized in \cite{visars} to the stack $\cH_{g}$ of hyperelliptic curves of genus~$g$, which has a presentation as a quotient $[X_{g}/\GL_{2}]$ (if $g$ is even), or $[X_{g}/(\gm\times\PGL_{2})]$ (if $g$ is odd), where $X$ is the space of smooth binary forms of degree~$2g+2$ in two variables; this allows to compute the Picard group of $\cH_{g}$. This description was applied in \cite{edidin-fulghesu-even} to compute the Chow ring of $\cH_{g}$ when $g$ is even.

In this paper we consider the stack $\T$ of smooth trigonal curves over a fixed base field $k$ of characteristic different from $2$ and $3$. Here by a trigonal curve over an algebraically closed field, we mean a smooth projective curve $C$ with a 3:1 map $C \arr P$, where $P$ is a curve isomorphic to $\PP^{1}$ (but the isomorphism is not part of the data). Consequently, we define $\T$ as follows. The objects of $\T$ over a $k$-scheme $S$ are of the form $(C \arr P \arr S)$, where $P \arr S$ is a smooth conic bundle (i.e., a family of $\PP^{1}$'s), $C \arr P$ a flat triple cover, and the composite $C \arr S$ is a smooth family of curves of genus $g$. There is a forgetful morphism $\T \arr \cM_{g}$, sending $(C \arr P \arr S)$ into the composite $C \arr S$; by Proposition~\ref{embed}, it is a locally closed embedding as soon as $g \geq 5$ (as might be expected, since uniqueness of the $g^{1}_{3}$ holds precisely in genus at least~$5$).

We give a description of $\T$ as a quotient $[U_{g}/\Gamma_{g}]$, where $\Gamma_{g}$ is a certain quotient of products of general linear groups and $X_{g}$ is an open subscheme of a $\Gamma_{g}$-equivariant vector bundle over an open subscheme of a representation of $\Gamma_{g}$ (see Theorem~\ref{quotient}). While the description above is considerably more complicated than that of $\cH_{g}$, we are able to exploit it to compute the Picard group of $\T$ when $g \neq 1$, thus proving the following.

\begin{theorem}\label{thm:picard}
Assume $g >1$. The Picard group of $\T$ is isomorphic to
   \[
   \begin{cases}
   \ZZ & \text{ if $g \not\equiv 0 \pmod 3$}\\
   \ZZ \oplus \ZZ/3\ZZ
      &\text{ if $g \equiv 0 \pmod 3$ and $g \not\equiv 3 \pmod 9$}\\
   \ZZ \oplus \ZZ/9\ZZ
      &\text{ if $g \equiv 3 \pmod 9$.}\\   
   \end{cases}
   \]

\end{theorem}

The calculation of the Picard group is only one of the possible applications of our explicit presentation, which we consider to be our main result. Another is the calculation of the rational Chow ring of $\cT_{g}$, which seems difficult, but feasible. We hope that this will be the subject of a subsequent paper.

It would be very interesting to extend this description to the closure of $\cT_{g}$ in $\overline{\cM}_{g}$; unfortunately, this does not seem to be easy. Even the case of $\cM_{2}$, treated in \cite{vistoli98}, presents problems: the obvious approach, using the Kirwan compactification, does not yield the stack $\overline{\cM}_{2}$ (although it probably gives the correct moduli space). We remark that the rational Picard group of the closure was computed in \cite{stankova}, by completely different methods.

Our approach has some similarities with that of \cite{sb-trigonal}.

The essential tools that we use are the following.

\begin{itemize}

\item A result of Miranda that describes a flat finite triple cover of a scheme as given by a locally free sheaf $E$ of rank two, with a section of $\sym^{3}(E) \otimes \det(E)^{\vee}$ (\cite{mir3}).

\item A description of the stack of globally generated locally free sheaves of fixed rank and degree on a projective line as a quotient stack (see Section~\ref{sec:glgen-conic}).

\end{itemize}

The ideas of this paper should extend to give a description of the stack of $d$-gonal curves, when one has a description of covers of degree $d$, in terms of ``generic'' linear algebra data. For $d = 4$ this is provided in \cite{caseke} and \cite{miranda-hahn-quadruple}, for $d = 5$ in \cite{casnati-degree5}. This will also be the subject of further work.

\subsection*{Description of contents}

Section~\ref{sec:definitions} contains the definition of the stack $\T$ that we are interested in, together with the proofs of two of its properties.

In Section~\ref{seclinalg} we state the results of Miranda and Casnati--Ekedahl in the form in which we will use them.

In Section~\ref{sec:glgen-conic} we give a description of the stack of globally generated locally free sheaves on a conic as a quotient stack. This might be of independent interest.

In Section~\ref{sectrigo} we put together the results of Miranda and those in Section~\ref{sec:glgen-conic} to obtain a description of a stack of triple covers $\OT$ that contains $\T$ as an open substack.

Section~\ref{secdet} is dedicated to the most intricate technical problem that we encounter, that of describing and analyzing the locus of singular curves in $\OT$. 

We conclude by using techniques of equivariant intersection theory together with all the results of the previous Sections to compute the Picard group of $\T$. Parts of the calculations have been made using the program \emph{Maple}, by Maplesoft.

\subsection*{Acknowledgments} We would like to thank the referees for very useful comments.

\section{Definitions and first results}\label{sec:definitions}

Throughout this paper we fix a base field $k$ of characteristic different from $2$ and $3$. Many of the things we do will work over more general rings, but for ease of exposition we disregard this fact.

For any non-negative integer $g$, we let $\Tbig$ be the category fibered in groupoids over  $\spec k$, whose sections over a scheme $S$ are of the form $(C \xarr{f} P \to S)$, where $P \arr S$ is a Brauer--Severi scheme of relative dimension~1, or, equivalently, a flat proper finitely presented morphism of schemes whose geometric fibers are isomorphic to $\PP^{1}$ (a \emph{conic bundle}, for short), $f\colon C \arr P$ is a finite flat finitely presented morphism of constant degree~$3$, such that the fibers of the composite $C \arr S$ have arithmetic genus~$g$. Conceptually, we should assume that $P$ and $C$ are algebraic space; but the invertible sheaves $\omega^{\vee}_{P/S}$ and $f^{*}\omega^{\vee}_{P/S}$ are ample over $S$, hence $P$ and $C$ are projective over $S$.

The fibered category $\Tbig$ is too large to be very useful; we are mostly interested in the full fibered subcategory $\T$ of objects $(C \xarr{f} P \to S)$ where the composition $(C \arr P \arr S)$ is supposed to be smooth. It is immediate to check that the embedding $\T \subseteq \Tbig$ is represented by open embeddings of schemes.

If $(C \xarr{f} P \to S)$ is in $\Tbig(S)$, the sheaf $f_{*}\cO_{C}$ is a locally free sheaf of rank~3 on $P$; the trace map $\tr\colon f_{*}\cO_{C} \arr \cO_{P}$ yields a splitting of the embedding $\cO_{P} \into f_{*}\cO_{C}$, because $3$ is invertible in $\cO(S)$; hence we can write $f_{*}\cO_{C} = \cO_{P} \oplus F$, where $F \into f_{*}\cO_{C}$ is the subsheaf of elements of trace~0, which is locally free of rank~2. As in \cite{mir3}, we call $F$ the \emph{\tsch module} of $(C \xarr{f} P \to S)$. We will mostly deal with its dual $E \eqdef F\dual$, the \emph{dual \tsch module}.

If $s\colon \spec \Omega \arr S$ is a geometric point, the geometric fiber $P_{s}$ of $P$ over $s$ is isomorphic to $\PP^{1}_{\Omega}$; hence the pullback of $E$ to $P_{s}$ will split as a direct sum of line bundles, which we denote by $\cO_{\PP^{1}_{\Omega}}(m) \oplus \cO_{\PP^{1}_{\Omega}}(n)$, with $m \leq n$. We call the pair of integers $(m,n)$ the \emph{splitting type} of $(C \xarr{f} P \to S)$ at the geometric point $s$. Since the splitting type only depends on the image of $s\colon \spec \Omega \arr S$, we can also talk of the splitting type of $(C \xarr{f} P \to S)$ at a point of $S$.

\begin{remark}
The in the case of smooth, or more generally Gorenstein, curves, the data of the splitting type is equivalent to the data of the genus and the Maroni invariant, which equals $m-2$ \cite{maroni, rosa-storh}.
\end{remark}

\begin{proposition}\call{prop:castelnuovo}
Let $(C \xarr{f} P \to S) \in \Tbig(S)$, $s\colon \spec \Omega \arr S$ be a geometric point and $(m,n)$ its splitting type at $s$. Then

\begin{enumeratea}

\itemref{1} $m + n = g + 2$, and

\itemref{2} if the geometric fiber $C_{s}$ is integral, then
   \[
   m, n \geq \frac{g+2}{3}.
   \]

\end{enumeratea}
\end{proposition}

\begin{proof}
For both statements we may base change to $\spec\Omega$, and assume that $f\colon C \arr \PP^{1}$ is a finite flat map of constant degree~$3$ defined over an algebraically closed field.

The equality in \refpart{prop:castelnuovo}{1} is immediate, by computing the Euler characteristic of $f_{*}{\cO}_{C}$ over $\PP^{1}$.

For the proof of \refpart{prop:castelnuovo}{2}, notice that because of part~\refpart{prop:castelnuovo}{1}, the inequalities to be proved are equivalent to $2m \geq n$ and $2n \geq m$; and for these see \cite[p.~1126]{mir3}.
\end{proof}

Call $\cM_{g}$ the stack of smooth curves of genus $g$ over $\spec k$. There is an obvious base-preserving functor $\T \arr \cM_{g}$ that sends $(C \arr P \to S)$ into $C \arr S$.

\begin{proposition}\label{embed}
If $g \geq 5$, then the functor $\T \arr \cM_{g}$ is a locally closed embedding.
\end{proposition}

This fails for $g \leq 4$, because for these values of $g$ there are smooth curves with more than one $g^{1}_{3}$, so the functor is not injective on isomorphism classes of geometric points.

\begin{proof}

Consider the open substack $\cM^{0}_{g}$ of $\cM_{g}$ consisting of non-hyperelliptic curves. Since a curve of genus at least~$3$ can not be both trigonal and hyperelliptic, the morphism $\T \arr \cM_{g}$ factors through $\cM^{0}_{g}$; we will prove that $\T \arr \cM^{0}_{g}$ is a closed embedding. For this it is enough to show that it is representable, proper, injective on geometric points and unramified.

To check that it is representable it is enough to check that, if $\Omega$ is an algebraically closed field and $(C \to P \to \spec\Omega)$ is a object of $\T$, the induced homomorphism from the automorphism group scheme of $(C \to P \to \spec\Omega)$ to the automorphism group scheme of $(C \to \spec\Omega)$ is injective. This follows immediately from the fact that the morphism $C \arr P$ is surjective.

The fact that this morphism is injective on geometric points is equivalent to the uniqueness of the $g^{1}_{3}$ for a smooth curve of genus at least~$5$. This is well known: if there are two morphisms $C \arr \PP^{1}$ of degree~$3$ which are different modulo the action of $\PGL_{2}$, the image $C'$ of the induced morphism $C \arr \PP^{1} \times \PP^{1}$ is birational to $C$ and bidegree $(3,3)$; hence the arithmetic genus of $C'$ is $4$ and $C$ is the normalization of $C'$, which implies $g \leq 4$.

Let us show that $\cT_g \arr \cM_{g}$ is unramified. This is equivalent to proving the following. If $(C \to P \to \spec A)$ and $(C \to P' \to \spec A)$ are two objects of $\cT_g(\spec A)$, where $A$ is an artinian ring with algebraically closed residue field, such that the two composites $C \arr \spec A$ coincide, then there exists an isomorphism $(C \to P \to \spec A) \simeq (C \to P' \to \spec A)$ in $\cT_g(\spec A)$ inducing the identity on $C$.

Call $\Omega$ the residue field of $A$. Since $P \simeq P' \simeq \PP^{1}_{A}$, the statement is equivalent to the following: given any two morphism $C \arr \PP^{1}_{A}$ of $A$-schemes which are finite of degree~$3$, these differ by an element of $\PGL_{2}(A)$. We know that the statement holds for $A = \Omega$, because of the previous point: hence we may assume that the restriction of the two morphisms $C \arr \spec A$ to $\spec \Omega$ coincide. By standard deformation theory, this can be rephrased as follows.

Let $f\colon C \arr \PP^{1}$ be a finite morphism of degree~$3$, where $C$ is a smooth curve of genus~$g$ over an algebraically closed field. Consider the tangent complex
   \[
   \bT_{C \to \PP^{1}} \eqdef
   \cdots \arr 0 \arr \rT_{C} \arr f^{*}\rT_{\PP^{1}} \arr 0 \arr \cdots,
   \]
where $f^{*}\rT_{\PP^{1}}$ is placed in degree~$0$. The tangent space to the deformation functor of the map $C \arr \PP^{1}$ is the $0\th$ hypercohomology group $\HH^{0}(C, \bT_{C \to \PP^{1}})$. There is an exact sequence
   \[
   \H^{0}(C, f^{*}\rT_{\PP^{1}})
   \arr \HH^{0}(C, \bT_{C \to \PP^{1}})
   \arr \H^{1}(C, \rT_{C})
   \]
where the last group map represents the natural homomorphism from the tangent space to the deformation functor of the map $C \arr \PP^{1}$ to the tangent space of the deformation functor of $C$. The tangent space to $\PGL_{2}$ is $\H^{0}(\PP^{1}, \rT_{\PP^{1}})$; hence it is enough to show that the natural homomorphism $\H^{0}(\PP^{1}, \rT_{\PP^{1}}) \arr \H^{0}(C, f^{*}\rT_{\PP^{1}})$ is surjective.

This is equivalent to saying that the homomorphism
   \[
   \H^{0}(\PP^{1}, \rT_{\PP^{1}})
   \arr \H^{0}(\PP^{1}, f_{*}f^{*}\rT_{\PP^{1}})
   = \H^{0}\bigl(\PP^{1}, f_{*}\cO_{C}\otimes\cO_{\PP^{1}}(2)\bigr)
   \]
induced by the embedding $\cO_{\PP^{1}}(2) = \cO_{\PP^{1}}\otimes\cO_{\PP^{1}}(2) \subseteq f_{*}\cO_{C}\otimes\cO_{\PP^{1}}(2)$ is surjective. Because of the splitting $f_{*}\cO_{C} = \cO_{\PP^{1}} \oplus F$, this is equivalent to saying that $\H^{0}\bigl(\PP^{1}, F \otimes \cO_{\PP^{1}}(2)\bigr) = 0$. If we denote by $(m,n)$ the splitting type of $C \arr \PP^{1}$, this is equivalent to $m \geq 3$. This last inequality follows from Proposition~\refall{prop:castelnuovo}{2}.

Finally we need to check that $\T \arr \cM^{0}_{g}$ is proper. It is representable, injective on geometric points and unramified: hence it is categorically injective, and so it is separated. It follows from the description in Section~\ref{sectrigo} that $\T$ is of finite type over $\spec k$, hence it is enough to check that the valuative criterion is satisfied. This means the following. Let $R$ be a discrete valuation ring with fraction field $K$, and let $C \arr \spec R$ be an object of $\cM^{0}_{g}(\spec R)$. Suppose that we are given an object $(C_{K} \to P_{K} \to \spec K)$ of $\T(\spec K)$, where $C_{K}$ is the restriction of $C$ to $\spec K$; then, after passing to a ramified extension of $R$, denoted by a standard abuse of notation also by $R$, there is an object $(C \to P \to \spec R)$ of $\T(\spec R)$ whose restriction to $\spec K$ is isomorphic to $(C_{K} \arr P_{K} \arr \spec K)$.

After passing to such an extension, we may assume that $P_{K}$ is isomorphic to $\PP^{1}_{K}$. The morphism $C_{K} \arr \PP^{1}_{K}$ is defined by a line bundle $L_{K}$ on $C_{K}$ of degree~$3$, with $\dim_{K}\H^{0}(C_{K}, L_{K}) \geq 2$. Extend the line bundle $L_{K}$ to a line bundle $L$ on $C$ (this is possible, because $C$ is regular). If we denote by $k$ the residue field of $R$, and by $C_{k}$ the closed fiber of $C \arr \spec R$, we have $\dim_{k} \H^{0}(C_{k}, L_{k}) \geq 2$. Since $C_{k}$ is not hyperelliptic, we have $\dim_{k}\H^{0}(C_{k}, L_{k}) = 2$, and that $L_{k}$ is generated by global sections. By Grauert's semicontinuity theorem it follows that $\H^{0}(C, L)$ is a free $R$-module of rank~$3$ and $L$ is generated by global section. Hence $L$ defines a morphism $C \arr \PP^{1}_{R}$ extending $C_{K} \arr \PP^{1}_{K}$, which is the required object of $\T(\spec R)$.
\end{proof}

\section{Triple covers via linear algebra}\label{seclinalg}

Here we recall some results on triple covers due to Miranda (\cite{mir3}) and Casnati\ddash Ekedahl (\cite{caseke}); see also \cite{faesti}.

Let $S$ be a scheme (or an algebraic space) over $\spec k$; denote by $\trip(S)$ the category whose objects are flat finite finitely presented maps $f\colon X \arr S$ of constant degree~$3$, which we call simply \emph{triple covers}, with the arrows given by isomorphisms over $S$. A triple cover is \emph{Gorenstein} when its geometric fibers are Gorenstein; if $S$ itself is a locally noetherian Gorenstein scheme, then $f$ is Gorenstein if and only if $X$ is Gorenstein.

Let us introduce another category $\trip'(S)$, with a very different description: its objects are pairs $(E,\alpha)$, where $E$ is a locally free sheaf of $\cO_{S}$-modules of constant rank~$2$, while $\alpha\in\H^{0}(S,\symthree{E})$ is a section. An arrow in $\trip'(S)$ from $(E,\alpha)$ to $(E', \alpha')$ is an isomorphism $\phi \colon E \simeq E'$ of $\cO_{S}$-modules, such that
   \[
   \symthree\phi\colon \symthree E \arr \symthree{E'}
   \]
carries $\alpha$ into $\alpha'$.

We define a functor $\trip(S) \arr \trip'(S)$ as follows. Given a triple cover $f \colon X \arr S$, let $F \subseteq f_{*}\cO_{X}$ the subsheaf consisting of elements of trace~$0$; since $3$ is invertible in $\cO(S)$, the trace map $f_{*}\cO_{X} \arr \cO_{S}$ gives a splitting $f_{*}\cO_{X} = \cO_{S} \oplus F$. Consider the composition
   \[
   \sym^{2} F \subseteq \sym^{2}f_{*}\cO_{X} \arr f_{*}\cO_{X} \arr F,
   \]
where the first inclusion is induced by the embedding $F \subseteq f_{*}\cO_{X}$, the second is given by the structure of commutative sheaf of algebras on $f_{*}\cO_{X}$, and the third is the projection. By duality, because of the canonical isomorphism $F \simeq F\dual \otimes \det F$, we obtain a homomorphism
   \[
   \sym^{2} F \otimes F \arr \det F,
   \]
which is easily seen to factor through the canonical projection $\sym^{2} F \otimes F \arr \sym^{3} F$, yielding a homomorphism $\alpha\colon \sym^{3} F \arr \det F$. Set $E \eqdef F\dual$. Using the customary identifications
   \[
   (\sym^3{F})\dual \simeq \sym^{3}E
   \quad\text{and}\quad
   \det F \simeq (\det E)\dual
   \]
(the former depending on the hypothesis that $1/6 \in \cO(S)$), we can think of $\alpha$ as a section of $\symthree{E}$. This defines a functor $\trip(S) \arr \trip'(S)$ that sends $f\colon X \arr S$ into $(E, \alpha)$.

\begin{theorem}[Miranda, Casnati-Ekedahl]\call{thm:describe-triple}\label{casna}\hfil
\begin{enumeratea}

\itemref{1} The functor $\trip(S) \arr \trip'(S)$ defined above is an equivalence.

\itemref{2} A triple covering is Gorenstein if and only if the corresponding section $\alpha \in H^0(S, \symthree{E})$ is nowhere vanishing.

\itemref{3} Assume that $f\colon X \arr S$ is Gorenstein. Consider the projection    
   \[
   \pi\colon \PP(E) \arr S
   \]
and the canonical isomorphism $\sym^{3}E \simeq \pi_{*}\cO_{\PP(E)}(3)$. We can think of $\alpha$ as a nowhere vanishing section of $\cO_{\PP(E\dual)}(3)\otimes \pi^{*}\det E$. Its zero locus $Y \subseteq \PP(E)$ is a triple cover of $S$, and is canonically isomorphic to $f\colon X \arr S$.

\end{enumeratea}
\end{theorem}


\section{The stack of globally generated locally free sheaves\\on a conic}\label{sec:glgen-conic}

In this section we will study two algebraic stacks over $\spec k$.

\begin{definition} Fix two non-negative integers $r$ and $d$.

\begin{enumeratea}

\item The objects of the category $\vectp$ are pairs $(S, E)$, where $S$ is a $k$-scheme and $E$ is a locally free sheaf of constant rank $r$ on $\pone_{S}$, whose restriction to any fiber of $\pone_{S} \arr S$ is globally generated and of constant degree $d$.

An arrow $(f, \phi)$ from $(S', E')$ consists of a morphism of $k$-schemes $f \colon  S' \arr S$, and an isomorphism of $\cO_{\pone_{S'}}$-modules $\phi\colon E' \simeq (\id_{\PP^{1}}\times f)^{*}E$. The composition is defined in the obvious fashion.

\item The objects of the category $\vect$ are pairs $(P \to S, E)$, where $S$ is a $k$-scheme, $P \arr S$ is a conic bundle, $E$ is a locally free sheaf of constant rank $r$ on $P$, whose restriction to any  fiber of $P \arr S$ is globally generated and of constant degree $d$.

An arrow $(f, F, \phi)$ from $(P' \to S', E')$ to $(P \to S, E)$ consists of a cartesian diagram of morphism of $k$-schemes
   \[
   \xymatrix{
   P' \ar[r]^{F}\ar[d] &P \ar[d]\\
   S' \ar[r]^{f}       &S\hsmash{,}
   }
   \]
where the columns are the obvious projections, and an isomorphism of $\cO_{P'}$-modules $\phi\colon E' \simeq F^{*}E$. The composition is defined in the obvious fashion.

\end{enumeratea}
\end{definition}

There are evident forgetful functors from $\vectp$ and from $\vect$ to the category of schemes, which make them into categories fibered in groupoids over $\spec k$. 

There is also a morphism of fibered categories $\vectp \arr \vect$ which sends $(S,E)$ into $(S, \pone_{S}, E)$, and an arrow $(f, \phi)$ into $(f, \id_{\PP^{1}}\times f, \phi)$. This is easily seen to make $\vectp$ into a $\PGL_{2}$-torsor over $\vect$.

Let us denote by $\mm$ the affine space over $k$ of $(r+d)\times d$ matrices $(\ell_{ij})$, where each $\ell_{ij}$ is a form of degree~$1$ in two indeterminates. We think of it as a scheme; in other words, $\rM_{r,g}$ represents the functor from $k$-schemes to sets, whose values on a $k$-scheme $S$ is the set of matrices $(r+d)\times d$ matrices $(\ell_{ij})$, where each $(\ell_{ij})$ is an element of $\H^{0}\bigl(\PP^{1}_{S}, \cO(1)\bigr)$. We identify such a matrix $(\ell_{ij})$ with the associated sheaf homomorphism
   \[
   \cO_{\pone_{S}}(-1)^{d} \xarr{(\ell_{ij})} \cO_{\pone_{S}}^{r+d}\hsmash{.}
   \]
We denote by $\om$ the open subscheme of $\mm$ parametrizing matrices $(\ell_{ij})$ with the property that the matrix $\bigl(\ell_{ij}(p)\bigr)$ has rank $d$ for all points $p \in (\AA^{2}\setminus \{0\})_{S}$; in other words, we assume, equivalently, that the morphism above is universally injective, or that is it injective and locally split, or that its dual is surjective. The cokernel of the universal homomorphism
   \[
   \cO_{\pone_{\om}}(-1)^{d} \arr \cO_{\pone_{\om}}^{r+d}
   \]
is a locally free sheaf of rank~$r$ on $\om$, which we denote by $E_{r,d}$.

\begin{proposition}\label{prop:codimension-1}
The codimension of $\mm \setminus \om$ into $\mm$ is at least~$r$.
\end{proposition}

\begin{remark}
In fact, one can show that $\mm \setminus \om$ has pure codimension~$r$.
\end{remark}

\begin{proof}
Let $V_{r,d}\subseteq \pone_{\mm}$ be the degeneracy locus of the universal homomorphism $\cO_{\pone_{\mm}}(-1)^{d} \arr \cO_{\pone_{\mm}}^{r+d}$, where this fails to have rank $d$. Its codimension is well-know to be at most $(r+d-d+1)(d-d+1) = r+1$; we need to show that it is in fact equal to $r+1$, and then the thesis will follow from the fact that $\mm \setminus \om$ is the image of $V_{r,d}$ along the projection $\pone_{\mm} \arr \mm$, which has relative dimension~1.

To do this, we can pull back to $(\AA^{2}\setminus\{0\})\times{\mm}$; call $\widetilde{V}_{r,d}$ the inverse image of $V_{r,d}$ in $(\AA^{2}\setminus\{0\})\times{\mm}$. There is a natural morphism $(\AA^{2}\setminus\{0\})\times{\mm} \arr \AA^{(r+d)d}$ into the affine space $\AA^{(r+d)d}$ of $(r+d)\times d$ matrices, obtained by sending a point $\bigl(p, (\ell_{ij})\bigr)$ of $(\AA^{2}\setminus\{0\})\times{\mm}$ into its evaluation $\bigl(\ell_{ij}(p)\bigr)$; this map is easily seen to be surjective, with fibers of constant dimension. Hence it is flat; but $\widetilde{V}_{r,d}$ is the inverse image of the locus in $\AA^{(r+d)d}$ of matrices with non-maximal rank, which is well known to have codimension~$r+1$.
\end{proof}

Let us denote by $\GL_{n}$ the group scheme $\GL_{n, k}$. There are three group schemes acting naturally on $\rM_{r,d}$, leaving $\om$ invariant.

The groups $\GL_{r+d}$ and $\GL_{d}$ act on $\mm$, the first by multiplication on the right and the second by multiplication on the left. There is also a left action of $\GL_{2}$ on $\H^{0}\bigl(\PP^{1}, \cO(1)\bigr) \simeq k^{2}$, which coincides with the tautological action. These three actions commute, hence induce an action of $\GL_{r+d} \times \GL_{d} \times \GL_{2}$ on $\mm$ and on $\om$, by the formula
   \begin{equation}\label{eq:action}
   (A,B,C)\cdot L = ACLB^{-1},\ \ L\in M_{r,d}.
   \end{equation}
There is an embedding $\gm \subseteq \GL_{r+d} \times \GL_{d} \times \GL_{2}$ by the formula
   \[
   t \arrto (\rI_{r+d},t\rI_{d},t^{-1}\rI_{2});
   \]
we denote by $T$ the image of this embedding; it is a central group subscheme of $\GL_{r+d} \times \GL_{d} \times \GL_{2}$. It is immediate to see that $T$ acts trivially on $\rM_{r.d}$, hence on $\om$; so the action of $\GL_{r+d} \times \GL_{d} \times \GL_{2}$ induces an action of the quotient
   \[
   \grp \eqdef (\GL_{r+d} \times \GL_{d} \times \GL_{2})/T
   \]
on $\rM_{r.d}$ and on $\om$. The action of $\GL_{r+d} \times \GL_{d}$ on $\mm$ is the restriction of the action of $\grp$.

Notice that there is a natural exact sequence of group schemes
   \[
   1 \larr \GL_{r+d} \times \GL_{d} \larr \grp \larr \PGL_{2} \larr 1,
   \]
where the homomorphism $\grp \arr \PGL_{2}$ is induced by the third projection of $\GL_{r+d} \times \GL_{d} \times \GL_{2}$ onto $\GL_{2}$.

Here is the main result of this section.

\begin{theorem}\label{thm:main-vector-bundles}
There are isomorphisms of fibered categories over $\spec k$
   \[
   \vect \simeq [\om/\grp]
   \]
and
   \[
   \vectp \simeq [\om/\GL_{r+d} \times \GL_{d}].
   \]
Under this isomorphisms, the natural morphism $\vectp \arr \vect$ corresponds to the morphism $[\om/\GL_{r+d} \times \GL_{d}] \arr [\om/\grp]$ induced by the embedding of $\GL_{r+d} \times \GL_{d}$ into $\grp$.
\end{theorem}

\begin{corollary}
The fibered categories $\vect$ and $\vectp$ are smooth geometrically integral algebraic stacks of finite type over $\spec k$ of dimensions $-r^{2}-3$ and $-r^{2}$ respectively.
\end{corollary}

This is not hard to show directly.

\begin{remark}
One can give a similar description of the stack of globally generated coherent sheaves on $\pone$ and on a conic; we simply need to substitute $\om$ with the larger open subset of $\mm$ of matrices $(\ell_{ij})$ which have rank~$d$ at the generic point of $\pone$.
\end{remark}

\begin{remark}
When $r = 1$ it is easy to see that the $\vectp[1,d]$ are all isomorphic to the classifying stack $\cB_{k}\gm$, while the $\vect[1,d]$ are isomorphic to $\cB_{k}(\gm \times \PGL_{2})$ when $d$ is even and to $\cB_{k}\GL_{2}$ when $d$ is odd; so in this case the presentations above are far from optimal. It is not clear to us if in the higher rank case the presentations above can be improved.
\end{remark}

\begin{proof}[Proof of Theorem~\ref{thm:main-vector-bundles}]
We will only prove the result for $\vect$; the part concerning $\vectp$ is similar but simpler, and is left to the reader.

We begin by giving an alternate description of $\vect$. Define the fibered category $\vectalt$ over $\spec k$, whose objects are of the type $(P \arr S, F_{0}, F_{1}, \phi)$, where
\begin{enumeratea}

\item $P \arr S$ is a conic bundle;

\item $F_{0}$ is a locally free sheaf of rank~$r+d$ on $P$ that is trivial on the geometric fibers of $P \arr S$;

\item $F_{1}$ is a locally free sheaf of rank~$d$ on $P$ that is isomorphic to $\cO(-1)^{d}$ on the geometric fibers of $P \arr S$;

\item $\phi\colon F_{1} \arr F_{0}$ is a locally split injection of sheaves of $\cO_{P}$-modules.

\end{enumeratea}

The arrows in $\vectalt$ are defined in the obvious way: an arrow from $(P' \arr S', F'_{0}, F'_{1}, \phi')$ to $(P \arr S, F_{0}, F_{1}, \phi)$ is a quadruple $(f, h, \theta_{0}, \theta_{1})$, where $f\colon S' \arr S$ and $h \colon P' \arr P$ are morphisms such that the diagram
   \[
   \xymatrix{
   P' \ar[r]^{h}\ar[d] & P\ar[d]\\
   S' \ar[r]^{f} & S
   }
   \]
is cartesian, and $\theta_{0}\colon F'_{0} \arr h^{*}F_{0}$ and $\theta_{1}\colon F'_{1} \arr h^{*}F_{1}$ are isomorphisms of $\cO_{P'}$-modules such that the diagram
   \[
   \xymatrix{
   F'_{1} \ar[r]^{\theta_{1}}\ar[d]^{\phi'} & h^{*}F_{1} \ar[d]^{h^{*}\phi}\\
   F'_{0} \ar[r]^{\theta_{0}} & h^{*}F_{0}
   }
   \]
commutes.

Let $(P\xarr{\pi} S, E)$ be an object of $\vect$. If $s\colon \spec\Omega \arr S$ is a geometric point, $P_{s} \simeq \PP^{1}_{\Omega}$ is the geometric fiber, and $E_{s}$ is the restriction of $E$ to $P_{s}$, we have $\dim_{\Omega}\H^{0}(P_{s}, E_{s}) = r+d$ and $\H^{1}(P_{s}, E_{s}) = 0$. Hence, by the standard base change theorems, the sheaf $\pi_{*}E$ is a locally free sheaf of rank $r+d$ on $S$, and its formation commutes with base change on $S$. The adjunction homomorphism $F_{0} \eqdef \pi^{*}\pi_{*}E \arr E$ is surjective; call $F_{1}$ its kernel, $\phi\colon F_{1} \arr F_{0}$ the embedding. Then $F_{1}$ is a locally free sheaf of rank~$d$ on $P$. If $s\colon \spec \Omega \arr S$ is a geometric point, we have an exact sequence of locally free sheaves on $P_{s}\simeq \PP^{1}_{\Omega}$
   \[
   0 \larr (F_{1})_{s} \stackrel{\phi_{s}}{\larr} (F_{0})_{s} \larr E_{s}
   \larr 0;
   \]
since the map $\H^{0}\bigl(P_{s}, (F_{0})_{s}\bigr) \arr \H^{0}\bigl(P_{s}, E_{s}\bigr)$ is an isomorphism and $\H^{1}\bigl(P_{s}, (F_{0})_{s}\bigr) = 0$, we have
   \[
   \H^{0}\bigl(P_{s}, (F_{1})_{s}\bigr)
   = \H^{1}\bigl(P_{s}, (F_{1})_{s}\bigr)
   = 0;
   \]
and this can only happen when $(F_{1})_{s} \simeq \cO(-1)^{d}$. Hence $(P, F_{0}, F_{1}, \phi)$ is an object of $\vectalt(S)$. This construction naturally extends to a base-preserving functor $\vect \arr \vectalt$ of categories fibered over $\spec k$.

We claim that this is an equivalence: the inverse functor $\vectalt \arr \vect$ is obtained by sending an object $(P \to S, F_{0}, F_{1}, \phi)$ of $\vectalt$ into the object $(P \to S, \coker \phi)$ of $\vect$. We leave the easy details to the reader.

Now we define another category fibered in groupoids over $\spec k$, which we call $\vectaux$. Its objects are triples $(P \arr S, F_{0}, F_{1})$, where

\begin{enumeratea}

\item $P \arr S$ is a conic bundle over a $k$-scheme $S$;

\item $F_{0}$ is a locally free sheaf of rank~$r+d$ on $P$ that is trivial on the geometric fibers of $P \arr S$;

\item $F_{1}$ is a locally free sheaf of rank~$d$ on $P$ that is isomorphic to $\cO(-1)^{d}$ on the geometric fibers of $P \arr S$.

\end{enumeratea}
These looks very much like $\vectalt$, except that the homomorphism $\phi$ is missing. The arrows are defined as quadruples $(f, h, \theta_{0}, \theta_{1})$ as for $\vectalt$, omitting the condition of the commutativity of the diagram involving the $\phi$'s.

Let us show that $\vectaux$ is a stack in the étale topology over $\spec k$. This follows from standard arguments. We leave it to the reader to check that it is a prestack, i.e., that arrows between pullbacks of two given objects form an étale sheaf. To check that descent data are effective, the only critical point is the effectiveness of descent data for conic bundles, since a in a conic bundle $P \arr S$ the space $P$ is supposed to be a scheme, and not simply an algebraic space: and this follows from descent for schemes with an ample invertible sheaf (\cite[Example~4.39]{fga-explained}).

Furthermore, any two objects of $\vectaux$ over the same scheme $S$ are locally isomorphic in the étale topology, and $\vectaux$ has a global object $\bigl(\PP^{1}, \cO_{\PP^{1}}^{r+d}, \cO_{\PP^{1}}(-1)^{d}\bigr)$ defined over $\spec k$; hence it is isomorphic to the stack of torsors of the sheaf in groups $\underaut_{k}\bigl(\PP^{1}, \cO_{\PP^{1}}^{r+d}, \cO_{\PP^{1}}(-1)^{d}\bigr)$ in the big étale site of $\spec k$. This sheaf sends each scheme $S$ into the group of automorphisms of $\bigl(\PP^{1}_{S}, \cO_{\PP^{1}_{S}}^{r+d}, \cO_{\PP^{1}_{S}}(-1)^{d}\bigr)$, whose objects are triples of the form $(\phi, f, g)$, where $\phi$ is an automorphism of $\PP^{1}_{S}$ (i.e., an element of $\PGL_{2}(S)$), $f$ is an isomorphism of $\cO_{\PP^{1}_{S}}^{r+d}$ with $\phi^*\cO_{\PP^{1}_{S}}^{r+d} = \cO_{\PP^{1}_{S}}^{r+d}$ (i.e., an element of $\GL_{r+d}\bigl(\cO(S)\bigr)$), and $g$ is an isomorphism of $\cO_{\PP^{1}_{S}}(-1)^{d}$ with $\phi^{*}\cO_{\PP^{1}_{S}}(-1)^{d}$. The composition is defined in the obvious way. There is a canonical homomorphism $\underaut_{k}\bigl(\PP^{1}, \cO_{\PP^{1}}^{r+d}, \cO_{\PP^{1}}(-1)^{d}\bigr) \arr \underaut_{k}(\PP^{1}) = \PGL_{2}$, whose kernel is canonically isomorphic to $\GL_{r+d} \times \GL_{d}$; the two factors act on $\cO_{\PP^{1}}^{r+d}$ and $\cO_{\PP^{1}}(-1)^{d}$ in the obvious way. 

We claim that $\underaut_{k}\bigl(\PP^{1}, \cO_{\PP^{1}}^{r+d}, \cO_{\PP^{1}}(-1)^{d}\bigr)$ is represented by $\grp$. In fact, the action of $\GL_{2} = \underaut_{k}\bigl(\PP^{1}, \cO_{\PP^{1}}(-1)\bigr)$ on $(\PP^{1}, \cO_{\PP^{1}}(-1))$ induces an action of $\GL_{2}$ on $(\PP^{1}, \cO_{\PP^{1}}^{r+d}, \cO_{\PP^{1}}(-1)^{d}\bigr)$, which commutes with the action of the subgroup $\GL_{r+d} \times \GL_{d}$ of $\underaut_{k}\bigl(\PP^{1}, \cO_{\PP^{1}}^{r+d}, \cO_{\PP^{1}}(-1)^{d}\bigr)$. This yields an action of $\GL_{r+d} \times \GL_{d} \times \GL_{2}$; it is easy to see that the subgroup $T \simeq \gm$ acts trivially. From this we obtain an action of $\Gamma_{r,d}$ on $\bigl(\PP^{1}, \cO_{\PP^{1}}^{r+d}, \cO_{\PP^{1}}(-1)^{d}\bigr)$, hence a homomorphism $\gamma\colon \Gamma_{r,d} \arr \underaut_{k}\bigl(\PP^{1}, \cO_{\PP^{1}}(-1)\bigr)$. We have a commutative diagram of non-abelian group schemes with exact rows
   \[
   \xymatrix{1 \ar[r]&
   \GL_{r+d} \times \GL_{d}\ar[r]\ar@{ = }[d]&
   \Gamma_{r,d} \ar[r]\ar[d]^{\gamma}&
   \PGL_{2}\ar[r]\ar@{ = }[d]& 1\\
   1 \ar[r]&
   \GL_{r+d} \times \GL_{d}\ar[r]&
   \underaut_{k}\bigl(\PP^{1}, \cO_{\PP^{1}}(-1)\bigr) \ar[r]&
   \PGL_{2}\ar[r]
   & 1
   }
   \]
which shows that $\gamma$ is an isomorphism. There is a tautological conic bundle $\cP_{r,d} \arr \vectaux$; this is the quotient $\pi\colon [\PP^{1}/\Gamma_{r,d}] \arr \cB_{k}\Gamma_{r,d} \simeq \vectaux$.

There are two tautological vector bundles $\cF_{0}$ and $\cF_{1}$ on $\cP_{r,d}$ of ranks $r+d$ and $d$, whose pullbacks along a morphism $S \arr \vectaux$ corresponding to an object $(P \arr S, F_{0}, F_{1})$ of $\vectaux(S)$ are $F_{0}$ and $F_{1}$ respectively. They correspond to the $\Gamma_{r,d}$-equivariant vector bundles $\cO_{\PP^{1}}^{r+d}$ and $\cO_{\PP^{1}}(-1)^{d}$ on $\PP^{1}$. By the standard base-change results, the formation of $\pi_{*}\underhom_{\cO_{\cP_{r,d}}}(\cF_{1}, \cF_{0})$ commutes with base change, so there is a vector bundle $\cH$ on $\cB_{k}\Gamma_{r,d}$ such that, given an object $(P \arr S, F_{0}, F_{1})$ of $\vectaux(S)$ corresponding to a morphism $S \arr \vectaux$, the sections of the pullback of $\cH$ to $S$ correspond to homomorphisms $F_{1} \arr F_{0}$. It is easy to see that the total space of the vector bundle $\cH$ is the quotient of $\hom_{\PP^{1}}(\cO_{\PP^{1}}(-1)^{d}, \cO_{\PP^{1}}^{r+d}) = \mm$ by the actions given by the formula~(\ref{eq:action}).

There is a morphism of fibered categories $\vectalt \arr \vectaux \simeq \cB_{k}\Gamma_{r,d}$ which sends $(P \arr S, F_{0}, F_{1}, \phi)$ into $(P \arr S, F_{0}, F_{1})$. This clearly factors through the total space of $\cH$, and it is immediate to see that it gives an isomorphism of $\vectalt$ with the open substack of the total space corresponding to the open subscheme $\om$ of $\mm$. This ends the proof of the theorem.
\end{proof}

\section{Vector bundles and trigonal curves}\label{sectrigo}

Now we connect stacks of trigonal curves with stacks of vector bundles on a conic. The stack $\Tbig$ is too large: we consider the open substack $\OT \subseteq \Tbig$ consisting of objects $(C \xarr{f} P \to S)$, whose splitting type $(m,n)$ at any point of $S$ satisfies $m, n \geq \frac{g+2}{3}$. Equivalently, we require that if $E$ is the dual \tsch module of $(C \xarr{f} P \to S)$, the sheaf $\symthree E$ is globally generated on any fiber of $C \arr S$. By the inequality of Proposition~\refall{prop:castelnuovo}{2} we see that $\T$ is an open substack of $\OT$.

Fix a splitting type $(m,n)$ with $n \geq m \geq 0$ and $m + n = g+2$. There is a locally closed substack $\cA_{m,n}$ of $\cV_{2, g+2}$, the full fibered subcategory whose objects are the object of $\cV_{2,g+2}$ which are étale-locally isomorphic to $\bigl(\PP^{1}_{S}, \cO(m)\oplus\cO(n)\bigr)$. 

\begin{lemma}\label{lem:codimension}
The locally closed substack $\cA_{m,n}$ of $\cV_{2, g+2}$ is smooth and irreducible. If $m < n$, its codimension is $n - m - 1$; if $m = n$, then it is an open substack.
\end{lemma}

\begin{proof}
Set $G \eqdef \underaut_{k}\bigl(\PP^{1}, \cO(m) \oplus \cO(n)\bigr)$. By descent theory, $\cA_{m,n}$ is isomorphic to the classifying stack $\cB_{k}G$. The obvious projection $G \arr \PGL_{2}$ is surjective, and its kernel is isomorphic to the sheaf of automorphisms of $\cO(m) \oplus \cO(n)$ as a sheaf on $\PP^{1}$. If $m = n$, this kernel is isomorphic to $\GL_{2}$, hence the dimension of $G$ is $7$, and the dimension of $\cA_{m,n}$ is $-7$, which equals the dimension of $\cV_{2,g+2}$ (as it must be, because $\cA_{m,n}$ is open in $\cV_{2,g+2}$). When $m < n$, the subsheaf $\cO(n)$ is preserved under any automorphism of $\cO(m) \oplus \cO(n)$. From this it easy to deduce that $G$ is a semi-direct product $(\gm\times\gm) \ltimes \hom_{\cO_{\PP^{1}}}\bigl(\cO(m), \cO(n)\bigr)$; hence its dimension is $n-m+3$. So the dimension of $G$ is $n - m + 6$, the dimension of $\cA_{m,n}$ is $-(n-m+6)$, and, since the dimension of $\cV_{2,g+2}$ is $-7$, its codimension is $n - m - 1$, as claimed.
\end{proof}

We will denote by $\widehat{\cV}_{g}$ the open substack of $\cV_{2, g+2}$ whose objects are pairs $(P \to S, E)$, where $P \arr S$ is a conic bundle and $E$ is a globally generated locally free sheaf of rank~$2$, degree~$g+2$, such that if $(m, n)$ is the splitting type of $E$ over a geometric fiber of $P \arr S$, then $m, n \geq \frac{g+2}{3}$.

From Lemma~\ref{lem:codimension} we deduce the following fact.

\begin{proposition}\label{prop:codimension-2}
Assume $g > 0$. Then the codimension of $\cV_{2, g+2} \setminus \widehat{\cV}_{g}$ in $\cV_{2, g+2}$ at least~$2$ everywhere.
\end{proposition}

We will shorten $\Gamma_{2,g+2}$ by $\Gamma_{g}$. Let $\widehat{\Omega}_{g} \subseteq \Omega_{2,g+2}$ be the inverse image of $\widehat{\cV}_{g} \subseteq \cV_{2,g+2}$ along the $\Gamma_{g}$-torsor $\Omega_{2,g+2} \arr \cV_{2,g+2}$. Assume that $g > 0$; then, by Propositions \ref{prop:codimension-1} and \ref{prop:codimension-2}, $\widehat{\Omega}_{g}$ is an open subscheme of the affine space $\rM_{2,g+2}$ whose complement has codimension at least~$2$. There is a universal conic bundle $\pi\colon \cP_{g} \arr \widehat{\cV}_{g}$, such that if $(P \to S, E)$ is an object of $\widehat{\cV}_{g}$, then $P \to S$ is the projection $S \times_{\widehat{\cV}_{g}} \cP_{g} \arr S$. There is also a tautological locally free sheaf $\cE_{g}$ on $\cP_{g}$ of rank~$2$. Since $\symthree {\cE_{g}}$ is globally generated along the fibers of $\cP_{g} \arr \widehat{\cV}_{g}$, the sheaf $\pi_{*}(\symthree {\cE_{g}})$ on $\widehat{\cV}_{g}$ is locally free of rank $2g+4$; furthermore, the formation of $\pi_{*}(\symthree {\cE_{g}})$ commutes with base change.

There is a natural projection $\OT \arr \widehat{\cV}_{g}$. By Theorem~\refall{thm:describe-triple}{1}, if $S$ is a scheme over $\spec k$, the category $\OT(S)$ is equivalent to the category of triples $(\pi\colon P \arr S, E, \alpha)$, where $(\pi\colon P \arr S, E)$ is an object of $\widehat{\cV}_{g}(S)$ and $\alpha$ is a section of $\symthree E$. This means that $\OT$ is the total space of the vector bundle on $\widehat{\cV}_{g}$ associated with $\pi_{*}(\symthree {\cE_{g}})$.

This description gives a presentation of $\OT$ as a quotient stack, starting from the isomorphism $\widehat{\cV}_{g} \simeq [\widehat{\Omega}_{g}/\Gamma_{g}]$. By definition, $\Gamma_{g}$ is a quotient of  $\GL_{g+4} \times \GL_{g+2} \times \GL_{2}$
which acts on $\pone$ via the projection onto $\GL_{2}$; this descends to an action of $\Gamma_{g}$ on $\pone$. This yields a natural action of $\Gamma_{g}$ on $\pone \times \widehat{\Omega}_{g}$; the quotient $[\pone \times \widehat{\Omega}_{g}/\Gamma_{g}]$ is $\cP_{g}$.

The tautological locally free sheaf $E_{g}$ is defined as the cokernel of the homomorphism
   \[
   \cO_{\pone_{\widehat{\Omega}_{g}}}(-1)^{g+2} \xarr{(\ell_{ij})}
   \cO_{\pone_{\widehat{\Omega}_{g}}}^{g+4},
   \]
where $(\ell_{ij}) \in \rM_{2, g+2}(\widehat{\Omega}_{g})$ is the tautological matrix of forms of degree~$1$. The sheaf $E_{g}$ is $\Gamma_{g}$-equivariant, and corresponds to the sheaf $\cE_{g}$ on $\cP_{g}$. The sheaf $\pi_{*}(\symthree{E_{g}})$ on $\widehat{\Omega}_{g}$ is locally free and its formation commutes with base change. Denote by $X_{g} \arr \widehat{\Omega}_{g}$ the total space of the vector bundle on $\widehat{\Omega}_{g}$ corresponding to $\pi_{*}(\symthree{E_{g}})$, that is, the relative spectrum of the symmetric algebra of $(\pi_{*}(\symthree{E_{g}}))\dual$. Note that $\Gamma_{g}$ acts naturally on $X_{g}$. Then our first main Theorem is a consequence of what we have just observed.

\begin{theorem}\label{quotient}
We have an equivalence of fibered categories $\OT\cong[X_{g}/\Gamma_{g}]$.

\end{theorem}

This allows us to compute the Picard group of $\OT$ for $g > 0$. This is the equivariant Picard group of the action of $\Gamma_{g}$ on $X_{g}$, which, by the homotopy invariance of the Chow groups, coincides with the equivariant Picard group of $\Gamma_{g}$ on $\widehat{\Omega}_{g}$. Since $\widehat{\Omega}_{g}$ is an invariant open subscheme of a representation of $\Gamma_{g}$ whose complement has codimension larger than $1$, this is the character group of $\Gamma_{g}$, which is isomorphic to $\ZZ \times \ZZ$.

But we are interested in the Picard group of the open subscheme $\T$. We will see in the next section that the complement $\OT \setminus \T$ is an irreducible hypersurface. The proof of Theorem~\ref{thm:picard} will be then concluded by computing its class in the character group of $\Gamma_{g}$.

\section{Detecting singular points in triple covers}\label{secdet}

\subsection{The stack of triple covers as a quotient stack.} Let us introduce some more notation. If $U$ is a finite free $k$-module, we still use $U$ to denote the associated scheme, that is, the spectrum of the $k$-algebra $\sym_{k}^{*}(U\dual)$. Thus, if $S$ is a $k$-scheme, the homomorphisms $S\arr U$ correspond to elements of $U \otimes_{k} \cO(S)$.

We will denote by $V$ the scheme associated with the $4$-dimensional $k$-vector space $\symthree[k]{k^{2}}$. The standard action of $\GL_{2}$ on $k^{2}$ induces an action of $\GL_{2}$ on $V$. If $F_{0}$ denotes the canonical rank~$2$ vector bundle on $\cB_{k}\GL_{2} \eqdef [\spec k/ \GL_{2}]$, the quotient stack $[V/\GL_{2}]$ is the total space of $\symthree{F_{0}}$; hence given a $k$-scheme $S$, the category $[V/\GL_{2}](S)$ is equivalent to the category of rank~$2$ locally free sheaves $E$ on $S$, together with a section of $\symthree{E}$, the arrows being isomorphisms of locally free sheaves preserving the sections.

We identify $k^{2}$ with $\H^{0}\bigl(\cO_{\PP^{2}}(1)\bigr)$; that is, we write the canonical basis as $x_{1}$, $x_{2}$ and think of an element $a_{1}x_{1} + a_{2}x_{2}$ of $k^{2}$ as a form $f(x)$ of degree~$1$ in two indeterminates, where $x=(x_1,x_2)$ is a row vector. (Actually, to be consistent we should think of $f$ as a column vector and write $f(x)$ as $xf$, but we will forgo this.) The action of $\GL_{2}$ on $k^{2}$ is then obtained by sending $\bigl(A,f(x)\bigr)$ into $f(xA)$. If we identify $\det k^{2}$ with $k$ by sending $x_{1} \wedge x_{2}$ to $1$, we can think of $\symthree{k^{2}}$ as the space $\H^{0}\bigl(\PP^{1}, \cO(3)\bigr)$ of forms $f(x)$ of degree~$3$ in two variables. The action of $\GL_{2}$ is given by
   \[
   \bigl(A,f(x)\bigr) \arrto \det(A)^{-1}f(xA).
   \]

Consider the stack in groupoids $\trip_{k} \arr \catsch{k}$ (see Section \ref{seclinalg}), such that for any $S \arr \spec k$ the objects of $\trip_{k}(S)$ consist of triple covers $X \arr S$. In view of Theorem \refall{thm:describe-triple}{1}, we see that $\trip_{k}$ is equivalent to $[V/\GL_{2}]$.

Now we need to treat the following problem. Given a smooth morphism $P \arr S$ with $1$-dimensional fibers and a morphism $P \arr [V/\GL_{2}]$, corresponding to a triple cover $X \arr P$, when is the composition $\rho\colon X \arr S$ smooth? Consider a geometric point $p\colon \spec\Omega \arr P$; the geometric fiber of $X$ over $p$ is of the form $\spec\Omega[x,y]/\bigl((x,y)^{2}\bigr)$ if $p$ maps to $0$ in $V$, while if it maps to a nonzero form $f(x) \in \H^{0}\bigl(\PP^{1}, \cO(3)\bigr)$, then the geometric fiber of $X$ over $p$ is the subscheme of $\PP^{1}_{\Omega}$ defined corresponding to the ideal $\bigl(f(x)\bigr) \subseteq \Omega[x_{1}, x_{2}]$, by Theorem \refall{thm:describe-triple}{3}. But knowing the fiber over $p$ does not determine whether there is a singular point of the map $X \arr S$ over $p$. For this we need to examine the inverse image in $X$ of the first order neighborhood of $p$ in $P$; we need a stack lying over $[V/\GL_{2}]$ that detects this information.

\subsection{The canonical thickening of a quotient stack}

If $R$ is a commutative $k$-algebra, we write $R[\epsilon]$ for the ring $R[\epsilon]/(\epsilon^{2})$.

Suppose that $X$ is a scheme over $k$. Set $X_{k[\epsilon]} \eqdef \spec k[{\epsilon}] \times_{\spec k} X$, and denote by $X_{\epsilon}$ the Weil transfer of $X_{k[\epsilon]}$ to $k$ via the embedding $k \subseteq k[\epsilon]$. The scheme $X_{\epsilon}$ represents the functor that sends a scheme $T \arr \spec k$ into the set of morphism of $k$-schemes $T_{k[\epsilon]} \arr X$. In other words, $X_{\epsilon}$ represents the functor $\underhom_{k}(\spec k[\epsilon], X)$. If $X$ is the scheme associated with a free $k$-module $W$, then $X_{\epsilon}$ is the scheme associated with the free $k$-module $W[\epsilon] \eqdef W \otimes k[\epsilon]$. 

Another description of $X_{\epsilon}$ is the following. Consider the sheaf $\Omega_{X/k}$ of K\"{a}hler differentials. Then $X_{\epsilon}$ is the relative spectrum of the symmetric algebra of $\Omega_{X/k}$ over $\cO_{X}$. Or, equivalently, $X_{\epsilon}$ is the normal cone of the diagonal embedding $X \into X\times_{k}X$.

There is a natural left action of $\gms{k}$ on $\spec k[\epsilon]$, induced by the left action on $k[\epsilon]$. If $R$ is a $k$-algebra, then $\gm(R)\eqdef R^{*}$ acts on the left on $R[\epsilon]$ as $u\cdot(a+\epsilon b) \eqdef a + \epsilon ub$; this induces a a left action of $\gms{R}$ on $X_\epsilon(R) = X(R[\epsilon])$, which is compatible with base change on $R$. In turn, this induces a left action of the group scheme $\gms{k}$ on $X_{\epsilon}$; this action corresponds to the grading of the symmetric algebra.

Now, suppose that $G \arr \spec k$ is a smooth finitely presented algebraic group acting on $X$. Denote by $\cF \eqdef [X/G]$ the quotient stack. There is a natural action of the group scheme $G_{\epsilon}$ over $X_{\epsilon}$: if $T$ is a scheme over $k$, then $G_{\epsilon}(T) \eqdef G(T_{k[\epsilon]})$ acts on $X_{\epsilon}(T) \eqdef X(T_{k[\epsilon]})$. Denote by $\frg$ the $k$-module corresponding to the normal bundle of the identity section $\spec k \arr G$; then by our general principle we identify $\frg$ with the scheme $\spec\sym_{k}\frg\dual$. We can think of $\frg$ as a group scheme via addition, on which $G$ acts with the adjoint representation. Then $G_{\epsilon}$ is canonically isomorphic to the semi-direct product $G\ltimes_{k} \frg$, because there is a canonical $G$-equivariant isomorphism $\Omega_{G/k} \simeq \cO_{G}\otimes_{k}\frg\dual$. The action of $\gms k$ on $G_{\epsilon}$ corresponds to the action of $\gm$ on $G\ltimes_{k} \frg$, leaving $G$ fixed and acting on $\frg$ by multiplication.

Call $H_G\eqdef \gms{k} \ltimes G_{\epsilon}$ the semi-direct product. This can also be thought of as the semi-direct product $(\gms{k} \times_{k} G) \ltimes \frg$, where $\gms{k}$ acts on $\frg$ by multiplication and $G$ by the adjoint action. There are two canonical projections $H_{G} \arr \gm$ and $H_{G} \arr G$.

The actions of $\gm$ and of $G_{\epsilon}$ on $X_{\epsilon}$ combine to give an action of $H_G$ on $X_{\epsilon}$. The quotient stack $[X_{\epsilon}/H_G]$ has a natural interpretation.

If $S$ is a scheme and $L$ a quasi-coherent sheaf on $S$, we denote by $S\generate{L}$ the relative spectrum over $S$ of the sheaf of algebras $\cO_{S} \oplus L$, where $L$ is a square-zero ideal. There is a canonical projection $S\generate{L} \arr S$, with a section $S \into S\generate{L}$ induced by the projection $\cO_{S} \oplus L \arr \cO_{S}$.

\begin{definition}\label{thick}
The \emph{canonical thickening of $\cF$} is the category $\cF^{(1)} \arr \catsch{k}$ fibered in groupoids defined as follows.  An object of $\cF^{(1)}$ is a triple $(S, L, \xi)$, where

\begin{enumeratea}

\item $S$ is a scheme over $k$,

\item $L$ is an invertible sheaf over $S$, and

\item $\xi$ is an object of $\cF(S\generate{L})$.

\end{enumeratea}

An arrow from $(S', L', \xi')$ to $(S, L, \xi)$ is a triple $(f, \phi, \psi)$, where

\begin{enumeratea}

\item $f\colon S' \arr S$ is a morphism of $k$-schemes,

\item $\phi\colon L' \simeq f^{*}L$ is an isomorphism of sheaves of $\cO_{S'}$-modules, and

\item $\psi\colon \xi' \arr \xi$ is an arrow of $\cF$ over the morphism $S'\generate{L'} \arr S\generate{L}$ induced by $f$ and $\phi$.

\end{enumeratea}

The composition is defined in the obvious way. The functor $\cF^{(1)} \arr \catsch{k}$ sends an object $(S, L, \xi)$ into $S$, and a morphism $(f, \phi, \psi)$ into $f$.

\end{definition}

It is straightforward to check that $\cF^{(1)}$ is fibered in groupoids.

There is a base-preserving functor $\cF^{(1)} \arr \cF$, sending an object $(S, L, \xi)$ into the restriction of $\xi$ to $S$ via the embedding $S \into S\generate{L}$ described above. Also, there is a tautological invertible sheaf $\cL$ on $\cF^{(1)}$, such that if $(S, L, \xi)$ is an object of $\cF^{(1)}(S)$, the pullback of $\cL$ along the corresponding morphism $S \arr \cF^{(1)}$ is isomorphic to $L$.

\begin{proposition}\label{isoquot}
The canonical thickening $\cF^{(1)}$ is canonically isomorphic to the quotient stack $[X_{\epsilon}/H_G]$.

Furthermore, the $\gm$-torsor on $\cF^{(1)}$ associated with the tautological $H_{G}$-torsor $X_{\epsilon} \arr \cF^{(1)}$ and the projection $H_{G} \arr \gm$ is the $\gm$-torsor corresponding to the tautological invertible sheaf on $\cF^{(1)}$; while the $G$-torsor coming from $X_{\epsilon} \arr \cF^{(1)}$ through the projection $H_{G} \arr G$ is the pullback of the tautological $G$-torsor $X \arr [X/G] =  \cF$ along the morphism $\cF^{(1)} \arr \cF$ above.
\end{proposition}

\begin{proof}

Let $L$ be an invertible sheaf on a scheme $S$, and denote by $\rho\colon L^{0} \arr S$ the associated $\gm$-torsor. A section of $L^{0}$ corresponds to a trivialization of $L$, which in turn yields an isomorphism of $S$-schemes $S_{k[\epsilon]} \simeq S\generate{L}$. This embeds $L^{0}$ into the scheme $\underhom_{S}(S_{k[\epsilon]}, S\generate{L})$; this last scheme is the scheme $S\generate{L}_{\epsilon}$, considered as an $S$-scheme.

Let us define a base-preserving functor $\cF^{(1)} \arr [X_{\epsilon}/H_G]$. Let $P\arr S\generate{L}$ be a $G$-torsor with a $G$-equivariant morphism $P \arr X$, giving an object of $\cF(S)$. Consider $P$ as an $S$-scheme through the composite $P \arr S\generate{L} \arr S$; call $Q$ the inverse image of $L^{0}$ in $P_{\epsilon}$, where the morphism $P_{\epsilon} \arr S\generate{L}_{\epsilon}$ is induced by the given morphism $P \arr S\generate{L}$. There are natural actions of $\gm$ and $G_{\epsilon}$ on $P_{\epsilon}$ described above; these leave $Q$ invariant, and are easily seen to induce an action of $H_G$. In this way we obtain a scheme $Q$ with an action of $H_G$ and an invariant morphism $Q \arr S$.

We claim that $Q$ is an $H_G$-torsor on $S$. The construction above commutes with base change on $S$, hence we may assume that $S = \spec k$, $L = \cO_{\spec k}$ and $P = G$. In this case we have $P_{\epsilon} = \underhom_{k}(\spec k[\epsilon], G_{k[\epsilon]}) = \underhom_{k}(\spec k[\epsilon], k[\epsilon])\times_{\spec k}G_{\epsilon}$, and $Q$ is isomorphic to $\gm \times_{\spec k} G_{\epsilon} = H_G$.

Now, the $G$-equivariant morphism $P \arr X$ induces an $H_G$-equivariant morphism $P_{\epsilon} \arr X_{\epsilon}$; when we restrict this to $Q$, the $H_G$-torsor $Q \arr S$ together with the resulting $H_G$-equivariant morphism $Q \arr X_{\epsilon}$ gives an object of $[X_{\epsilon}/H_G](S)$. This extends immediately to arrows, and gives a base-preserving functor $\cF^{(1)} \arr [X_{\epsilon}/H_G]$.

To go in the opposite direction, start from an $H_G$-torsor $Q \arr S$. The quotient $Q/G_{\epsilon} \arr S$ by the subgroup $G_{\epsilon} \subseteq H_G$ is a $\gm$-torsor, which is of the form $\rho\colon L^{0} \arr S$ for a canonically defined invertible sheaf $L$ on $S$. The pullback of $L$ to $L^{0}$ has a tautological section, which gives an isomorphism $\cO_{L^{0}} \simeq \rho^{*}L$, inducing an isomorphism of $S$-schemes $L^{0}\generate{\rho^{*}L} \simeq L^{0}_{k[\epsilon]}$. Then the projection $Q \arr L^{0}$ is $H_G$-equivariant, when we let $H_G$ act on $\gm$ via the projection $H_G \arr H_G/G_{\epsilon} = \gm$.

There are two actions of $\gm$ on $L^{0}_{k[\epsilon]}$; also, there are two natural actions of $\gm$ on $L^{0}\generate{\rho^{*}L} = L^{0}\times_{S}S\generate{L}$, one induced by the multiplication on $L^{0}$, the other on $\rho^{*}L$. The diagonal action on $L^{0}_{k[\epsilon]}$ corresponds to the action on the first component of $L^{0}\times_{S}S\generate{L}$; hence the quotient $L^{0}_{k[\epsilon]}$ by the diagonal action is $S\generate{L}$. We let $H_G$ act on $Q_{k[\epsilon]}$ as the product of the given action on $H_G$ and the action on $\spec k[\epsilon]$ through the projection $H_G \arr \gm$. Then the projection $Q_{k[\epsilon]} \arr L^{0}_{k[\epsilon]}$ is $H_G$-equivariant. We may think of $Q_{k[\epsilon]} \arr L^{0}_{k[\epsilon]}$ as a $\gm$-equivariant $G_{\epsilon}$-torsor, where $\gm$ acts on $G_{\epsilon}$ as described above. The quotient $Q_{k[\epsilon]}/\frg$ is a $\gm$-equivariant $G$-bundle on $L^{0}_{k[\epsilon]}/\gm = S\generate{L}$; hence, since the action of $\gm$ on $L^{0}_{k[\epsilon]}$ is free, we obtain a $G$-torsor $P \eqdef Q_{k[\epsilon]}/(\gm\ltimes \frg)$ over $L^{0}_{k[\epsilon]}/\gm = S\generate{L}$.

Now, the morphism $Q \arr X_{\epsilon}$ corresponds, by definition, to a morphism of $k$-schemes $Q_{k[\epsilon]} \arr X$. This is $H_G$-equivariant, when we let $H_G$ act on $X$ on via the natural projection $H_G \arr G$. The kernel of the homomorphism $H_G \arr G$, which is the semi-direct product $\gm \ltimes \frg$, acts trivially on $X$; hence there is an induced $G$-equivariant morphism $P = Q_{k[\epsilon]}/(\gm\ltimes \frg) \arr X$. The data of the $G$-torsor $P \arr S\generate{L}$ and the morphism $P \arr X$ give an object of $\cF^{(1)}$. This construction extends to the arrows in the obvious way, and it defines a base-preserving functor $[X_{\epsilon}] \arr \cF^{(1)}$.

We need to check that these two functors are inverse to each other. This is straightforward and left to the reader.

The last statement follows from the construction.
\end{proof}

\subsection{The thickened stack of triple covers}\label{subsec:thickened-triple} Now we apply the preceding construction to the case we are interested in. If $R$ is a $k$-algebra, then $\GLe(R) = \GL_{2}(R[\epsilon])$. Thus each element of $\GLe(R)$ will be written as $A + \epsilon B$, where $A \in \GL_{2}(R)$ and $B \in \rM_{2}(R)$. Let us write $\Hg$ for the semi-direct product $\gm \ltimes \GLe(k)$.

From the isomorphism $\trip_{k} \simeq [V/\GL_{2}]$, where $V = \symthree[k]{k^{2}}$, and from Proposition~\ref{isoquot} we obtain an equivalence
   \[
   \trip_{k}^{(1)} \simeq [\Ve/\Hg],
   \]
where $\Ve$ is the scheme associated with the free  $k$-module of rank~$8$
   \[
   \symthree[{k[\epsilon]}] {k[\epsilon]^{2}}.
   \]

We will need formulas for the action of $\Hg$ on $\Ve$. As before, we think of $k[\epsilon]^{2}$ as the free $k[\epsilon]$-module over the two indeterminates $x_{1}$ and $x_{2}$; so we write elements of $k[\epsilon]^{2}$ as forms $(f+\epsilon g)(x)$ of degree~$1$ in the indeterminates $x_{1}$, $x_{2}$ and coefficients in $k[\epsilon]$, and think of $x = (x_{1}, x_{2})$ as a row vector. Then $\GLe$ acts on $k[\epsilon]^{2}$ by the formula
   \begin{align*}
   (A + \epsilon B) \cdot (f + \epsilon g)(x) &=
      (f + \epsilon g)\bigl(x(A + \epsilon B)\bigr)\\
      &= f(xA) + \epsilon\bigl(f(xB) + g(xA)\bigr).
   \end{align*}

The $k[\epsilon]$-module $\det_{k[\epsilon]}k[\epsilon]^{2}$ is identified with $k[\epsilon]$ by sending $x_{1} \wedge x_{2}$ to $1$; hence $\Ve$ is identified with the space $\sym^{3}_{k[\epsilon]}k[\epsilon]^{2}$ of forms of degree~$3$ in $x_{1}$, $x_{2}$, which we write as $f(x) + \epsilon g(x)$, where $f$ and $g$ are forms of degree~$3$ with coefficients in $k$. The action of $\GLe$ on $\Ve$ is written as
   \begin{align*}
   (A + \epsilon B) \cdot\bigl(f(x) + \epsilon g(x)\bigr) &=
    \det(A + \epsilon B)^{-1}
    \bigl(f(xA + \epsilon xB) + \epsilon g(xA)\bigr)\\
   &= \det(A + \epsilon B)^{-1}\bigl(f(xA) + \epsilon\bigl(xB \rJ_{f}(xA)
      + g(xA)\bigr)\bigr),
   \end{align*}
where
   \[
   \rJ_{f}(x)\eqdef 
   \begin{pmatrix}
   \frac{\partial f(x)}{\partial x_{1}}\\\vspace*{-7pt}\\
   \frac{\partial f(x)}{\partial x_{2}}
   \end{pmatrix}
   \]
is the jacobian matrix of $f$.

The action of $\gm$ on $\Ve$ is given by the formula $ \alpha\cdot(f + \epsilon g) \eqdef f + \epsilon \alpha g$.

\subsection{The evaluation map $\Lambda\colon \cP_g' \arr \trip_k^{(1)}$.} Call $\cP_g'$ the pullback of the universal conic bundle $\cP_g \arr \widehat{\cV}_{g}$ along the natural projection $\OT\arr \widehat{\cV}_{g} \subseteq \cV_{2,g+2}$; there is a tautological morphism
   \[
   \Lambda_{0}\colon \cP'_{g} \arr \trip_{k}.
   \]
If $S$ is a scheme over $k$, elements $(\pi\colon P\rightarrow S,E,\alpha)$ of $\OT (S)$ determine morphisms $P\rightarrow \trip_k$. If $I$ is the sheaf of ideals of the diagonal embedding $P \into P\times_{S}P$, we have $I/I^{2} = \Omega_{P/S}^{1}$. We have two isomorphisms of the closed subscheme of $P\times_{S}P$ corresponding to $I^{2}$ with $P\generate{\Omega_{P/S}^{1}}$; we choose the one under which the natural morphism $P\generate{\Omega_{P/S}^{1}} \arr P$ is induced by the second projection $P\times_{S}P \arr P$. Consider the composite $P\generate{\Omega_{P/S}^{1}} \arr P \arr \trip_{k}$, where the first arrow is induced by the first projection $P\times_{S}P \arr P$; this gives a morphism $P \arr \trip_{k}^{(1)}$. This construction is obviously compatible with pullbacks on $S$, hence it induces a morphism
   \[
   \Lambda\colon \cP'_{g} \arr \trip_{k}^{(1)} \simeq [\Ve/\Hg]
   \]
lifting $\Lambda_{0}\colon \cP'_{g} \arr \trip_{k}$. 

The two projections $\Hg \arr \gm$ and $\Hg \arr \GL_{2}$ induce, via pullback along $\Lambda$, two locally free sheaves on $\cP'_{g}$; the first is the canonical bundle $\Omega^1_{\cP'_g/\OT}$, the second is the tautological rank two locally free sheaf $\cE'_{g}$ on $\cP'_{g}$.

\begin{proposition}\label{prop:valuation-map}\hfil
\begin{enumeratea}

\item If $g \not\equiv 1 \pmod 3$, the morphism $\cP'_{g} \arr [\Ve/\Hg]$ is smooth with irreducible geometric fibers.

\item If $g \equiv 1 \pmod 3$, then there exist two irreducible closed substacks $\cA$  and $\cB$ of $\cP'_{g}$, of codimension  $\frac{g-1}{3}$ and $\frac{g+20}{3}$ respectively, with $\cB \subseteq \cA$, such that the restriction $\cP'_{g} \setminus \cA \arr [\Ve/\Hg]$ of the morphism $\cP'_{g} \arr [\Ve/\Hg]$ is smooth with irreducible geometric fibers, while  $\cA \setminus \cB \arr [\Ve/\Hg]$ is flat.

\end{enumeratea}
\end{proposition}

From this one deduces the following.

\begin{corollary}\label{cor:integral-inverse-image}
Let $W$ be a geometrically integral $\Hg$-invariant subscheme of $\Ve$ of codimension~$c$. If $g \neq 1$, then the inverse image $\Lambda^{-1}[W/\Hg] \subseteq \cP'_{g}$ is geometrically integral of codimension~$c$, except when $g = 4$ and $W = \{0\}$.
\end{corollary}

\begin{proof}
If $g \not\equiv 1 \pmod 3$, the result follows immediately from the Proposition. If $g \equiv 1 \pmod 3$ and $g \neq 1$, then the intersection $\Lambda^{-1}[W/\Hg] \cap (\cP'_{g} \setminus \cA)$ is smooth and irreducible, so it is enough to show that $\Lambda^{-1}[W/\Hg]$ does not have any components contained in $\cA$. Now, the codimension of $\Lambda^{-1}[W/\Hg] \cap(\cA \setminus\cB)$ in $\cA \setminus\cB$ equals the codimension of $W$ in $V$, by the second part of the proposition. The only problem may occur when the codimension of $\cB$ in $\cP'_{g}$, which is $\frac{g+20}{3}$, is less than or equal to the codimension of $W$ in $V$, which is at most 8; and this can occur only if $g = 4$ and $W = \{0\}$.
\end{proof}

\begin{proof}[Proof of Proposition~\ref{prop:valuation-map}]
Denote by $\widehat{\cW}_{g}$ the open substack of $\widehat{\cV}_{g}$, defined by the condition that an object $(P \to S, E)$ of $\widehat{\cV}$ is in $\widehat{\cW}_{g}$ if and only if the splitting type $(m, n)$ of $E$ on each fiber of $P \to S$ is such that $m$, $n > \frac{g+2}{3}$. Call $\cQ_{g}$ and $\cQ'_{g}$ the inverse images of $\widehat{\cW}_{g}$ in $\cP_{g}$ and $\cP'_{g}$ respectively. If $g \not\equiv 1 \pmod 3$ then $\cQ'_{g} = \cP'_{g}$; while if $g \equiv 1 \pmod 3$ the codimension of the complement of $\cQ'_{g}$ in $\cP'_{g}$ is $\frac{g-1}{3}$, by Lemma~\ref{lem:codimension}. Hence the following Lemma implies the first statement.

\begin{lemma}
The restriction $\cQ'_{g} \arr [\Ve/\Hg]$ is smooth with irreducible geometric fibers.
\end{lemma}

\begin{proof}
By Proposition~\ref{isoquot} applied to the case $X = \spec k$, an object of the classifying stack $\cB_{k}\Hg$ can be described as given by a $k$-scheme $S$, an invertible sheaf $L$ on $S$, and a locally free sheaf of rank $2$ on $S\generate{L}$. We can produce a morphism $\cQ_{g} \arr \cB_{k}\Hg$ as follows. An object of $\cQ_{g}$ is given by a scheme $S$, a conic bundle $P \arr S$, a section $\sigma\colon S \arr P$, and a rank~$2$ locally free sheaf $F$ on $P$. Let $I$ be the sheaf of ideals of the closed embedding $\sigma\colon S \arr P$. The conormal bundle $I/I^{2}$ is $\sigma^{*}\Omega^{1}_{P/S}$; furthermore, the first order neighborhood $\spec \cO_{P}/I^{2}$ of $S$ inside $P$ is canonically isomorphic to $\spec(\cO_{P} \oplus \sigma^{*}\Omega^{1}_{P/S}) = S\generate{\sigma^{*}\Omega^{1}_{P/S}}$ (the projection $\cO_{P}/I^{2} \arr \cO_{P}$ is split by the morphism $P \arr S$). We associate with this object of $\cQ_{g}$ the object of $\cB_{k}\Hg$ given by the scheme $S$, the invertible sheaf $\sigma^{*}\Omega^{1}_{P/S}$ and the restriction of $F$ to $S\generate{\sigma^{*}\Omega^{1}_{P/S}}$ via the embedding $S\generate{\Omega^{1}_{P/S}} \into P$.

This morphism $\cQ_{g} \arr \cB_{k}\Hg$ is smooth, with irreducible geometric fibers. To check this, we can base change to $\spec k$ via the usual morphism $\spec k \arr \cB_{k}\Hg$ given by the trivial $\Hg$-torsor on $\spec k$; the morphism $\cQ_{g} \arr \spec k$ is smooth with irreducible geometric fibers, and the fibered product $\cQ_{g}\times_{\cB_{k}\Hg}\spec k$ is an $\Hg$-torsor over $\cQ_{g}$. Since the group $\Hg$ is smooth and connected, the composite $\cQ_{g}\times_{\cB_{k}\Hg}\spec k \arr \cQ_{g} \arr \spec k$ is smooth with irreducible geometric fibers, as claimed.

Thus we have morphisms $\cQ'_{g} \arr [\Ve/\Hg]$ and $\cQ_{g} \arr \cB_{k}\Hg$. It is easy to see that the two composites $\cQ'_{g} \arr \cQ_{g} \arr \cB_{k}\Hg$ and $\cQ' \arr [\Ve/\Hg] \arr \cB_{k}\Hg$ are canonically isomorphic; in other words, we have a commutative diagram
   \[
   \xymatrix{
   \cQ'_{g} \ar[r]\ar[d]& [\Ve/\Hg] \ar[d]\\
   \cQ_{g}\ar[r] & \cB_{k}\Hg.
   }
   \]
The resulting morphism $\cQ'_{g} \arr \cQ_{g}\times_{\cB_{k}\Hg}[\Ve/\Hg]$ can be described as follows. An object $(P\to S, \sigma\colon S \to P, F,\alpha)$ of $\cQ_{g}\times_{\cB_{k}\Hg}[\Ve/\Hg]$ corresponds to the following data:

\begin{enumeratea}

\item a conic bundle $P \arr S$ with a section $\sigma\colon S \arr P$,

\item a locally free $F$ of rank~$2$ on $P$, and

\item an element $\alpha$ of $\H^{0}\bigl(S\generate{\sigma^{*}\Omega^{1}_{P/S}}, \symthree{F} \mid_{S\generate{\sigma^{*}\Omega^{1}_{P/S}}}\bigr)$, where $S\generate{\sigma^{*}\Omega^{1}_{P/S}}$ is considered as embedded into $P$ in the way described above.

\end{enumeratea}

An object $(P\to S, \sigma\colon S \to P, F, \beta)$ of $\cQ'_{g}$ consists of the following:

\begin{enumeratea}

\item a conic bundle $P \arr S$ with a section $\sigma\colon S \arr P$,

\item a locally free $F$ of rank~$2$ on $P$, and

\item an element $\beta$ of $\H^{0}(P, \symthree{F})$.

\end{enumeratea}

The morphism $\cQ'_{g} \arr \cQ_{g}\times_{\cB_{k}\Hg}[\Ve/\Hg]$ associates with an object $(P\to S, \sigma\colon S \to P, F, \beta)$ the element $(P\to S, \sigma\colon S \to P, F, \beta\mid_{S\generate{\sigma^{*}\Omega^{1}_{P/S}}})$. From this description it is clear that both $\cQ'_{g}$ and $\cQ_{g}\times_{\cB_{k}\Hg}[\Ve/\Hg]$ are vector bundles on $\cQ_{g}$, and the morphism $\cQ'_{g} \arr \cQ_{g}\times_{\cB_{k}\Hg}[\Ve/\Hg]$ is linear. We claim that it also surjective.

It is enough to check this it on geometric points. So, let us $S = \spec \Omega$ be the spectrum of an algebraically closed field, and let $(P\to S, \sigma\colon S \arr P, F, \alpha)$ be a point of $\cQ_{g}$ over $S$. Then $P \simeq \PP^{1}_{\Omega}$; denote by $p \in P$ the image of $\sigma$, so that $S\generate{\sigma^{*}\Omega^{1}_{P/S}}$ is the divisor $2p$. Furthermore, $F \simeq \cO(m) \oplus \cO(n)$, with $m$, $n > \frac{g+2}{3}$, $m+n = g+2$. Then $\symthree{F}$ decomposes as $\cO(2m-n)\oplus \cO(m) \oplus \cO(n) \oplus \cO(2n-m)$, and $2n-m \geq n \geq m \geq 2m-n >0$; so $\H^{1}\bigl(P, \symthree{F}(-2p)\bigr) = 0$, and the restriction homomorphism
   \[
   \H^{0}(P, \symthree{F}) \arr \H^{0}(2p, \symthree{F}\mid_{2p})
   \]
is surjective, which proves what we want.

So the morphism $\cQ'_{g} \arr \cQ_{g}\times_{\cB_{k}\Hg}[\Ve/\Hg]$ is a surjective homomorphism of vector bundles over $\cQ_{g}$, hence it is smooth with irreducible geometric fibers. The projection $\cQ_{g}\times_{\cB_{k}\Hg}[\Ve/\Hg] \arr [\Ve/\Hg]$ is obtained by base change from $\cQ_{g} \arr \cB_{k}\Hg$, which is smooth with irreducible geometric fibers; hence it is smooth with irreducible geometric fibers. So the composite
   \[
   \cQ'_{g} \arr \cQ_{g}\times_{\cB_{k}\Hg}[\Ve/\Hg]
      \arr [\Ve/\Hg]
   \]
is irreducible with geometric fibers, as claimed.
\end{proof}

For the second statement, let $\cA'$ be the complement of $\widehat{\cW}_{g}$ in $\widehat{\cV}_{g}$: this is the integral closed substack of $\widehat{\cV}_{g}$ denoted by $\cA_{\frac{g+2}{3}, \frac{2g+4}{3}}$ in the proof of Lemma~\ref{lem:codimension}. Set $(m, n) \eqdef \bigl(\frac{g+2}{3}, \frac{2g+4}{3}\bigr)$. From that proof we see that $\cA'$ is isomorphic to $\cB_{k}G$, where $G \eqdef \underaut_{\spec k}\bigl(\PP^{1}, \cO(m)\oplus \cO(n))$. Denote by $\cA$ the inverse image of $\cA'$ in $\cP'_{g}$; it is easy to see that $\cA$ is the quotient
   \[
   [\H^{0}\bigl(\PP^{1},
      \symthree{\bigl(\cO(m) \oplus \cO(n)\bigr)}\bigr)/G'],
   \]
where $G'$ is the stabilizer of the point $(1:0)$ of $\PP^{1}$ under the action of $G$ on $\PP^{1}$. 

Set
   \begin{align*}
   E &\eqdef \H^{0}\bigl(\PP^{1},
      \symthree{\bigl(\cO(m) \oplus \cO(n)\bigr)}\bigr)\\
   &= \H^{0}\bigl(\PP^{1}, \cO(3m) \oplus \cO(2m) \oplus \cO(m)
      \oplus \cO\bigr).
   \end{align*}

The projection $E \arr \cA$ is a $G'$-torsor, and $G'$ is smooth and connected: hence it is enough to show that the composite
   \[
   E \arr \cA \subseteq \cP'_{g}\arr [\Ve/\Hg],
   \]
which we denote by $\Phi\colon E \arr [\Ve/\Hg]$, is smooth with irreducible geometric fibers outside of a $G'$-invariant irreducible subscheme $T$ of $E$ of codimension~$7$. In this way we can set $\cB \eqdef [T/G']  \subseteq \cA$, and the codimension of $\cB$ in $\cP'_{g}$ is $\frac{g-1}{3} + 7 = \frac{g+20}{3}$.

Let us us denote by $\{0\}$ the origin in $\Ve$, considered as a closed subscheme. Consider the closed substack $[\{0\}/\Hg] \subseteq [\Ve/\Hg]$. The following Lemma completes the proof of the Proposition.
\end{proof}

\begin{lemma}
The inverse image $\Phi^{-1}\bigl([\{0\}/\Hg]\bigr)$ is irreducible of codimension~$7$ in $E$, and the restriction
   \[
   E \setminus \Phi^{-1}\bigl([\{0\}/\Hg]\bigr) \arr [\Ve \setminus \{0\}/\Hg]
   \]
is flat.
\end{lemma}

\begin{proof}
The morphism $\Phi\colon E \arr [\Ve/\Hg]$ factors through $\Ve$, as follows.

Let us use $t_{0}$ and $t_{1}$ for the homogeneous coordinates on $\PP^{1}$, as set $t \eqdef t_{1}/t_{0}$. Let us write an element of $E$ in the form $(\phi_{0}, \phi_{1}, \phi_{2}, \phi_{3})$, where the $\phi_{i}$'s are polynomials in $t$ of degrees $3m$, $2m$, $m$ and $0$ respectively. Define a linear map $E \arr \Ve$ by the formula
   \begin{align*}
   (\phi_{0}, \phi_{1}, \phi_{2}, \phi_{3}) &\arrto
   \phi_{0}(\epsilon)x_{1}^{3} + \phi_{1}(\epsilon)x_{1}^{2}x_{2}
      +\phi_{2}(\epsilon)x_{1}x_{2}^{2}
      +\phi_{3}(\epsilon)x_{2}^{3}\\
   &\quad = f + \epsilon g,
   \end{align*}
where we have set
   \begin{align*}
   f \eqdef \phi_{0}(0)x_{1}^{3} + \phi_{1}(0)x_{1}^{2}x_{2}
      +\phi_{2}(0)x_{1}x_{2}^{2}
      +\phi_{3}(0)x_{2}^{3}\\
\intertext{and}
   g \eqdef \phi_{0}'(0)x_{1}^{3} + \phi_{1}'(0)x_{1}^{2}x_{2}
      +\phi_{2}'(0)x_{1}x_{2}^{2}
      +\phi_{3}'(0)x_{2}^{3}.
   \end{align*}
The image of $E$ into $V_{\epsilon}$ is easily seen to be the hyperplane
   \[
   U \eqdef \{f+\epsilon g \in \Ve\mid g(0,1) = 0\}.
   \]
The homomorphism $E \arr U$ is linear and surjective, hence it is smooth with irreducible fibers. The inverse image of $[\{0\}/\Hg]$ in $U$ is the origin $\{0\}$, which is irreducible of codimension~$7$: this proves the first part of the Lemma.

For the second part, it is enough to show that the restriction
   \[
   U\setminus\{0\} \subseteq \Ve\arr [\Ve\setminus\{0\}/\Hg]
   \]
is flat. Set $U^{0} \eqdef U\setminus\{0\}$, $\Ve^{0} \eqdef \Ve \setminus \{0\}$. We have a cartesian diagram
   \[
   \xymatrix{
   \Hg \times U^{0} \ar[r]\ar[d]& \Ve^{0} \ar[d]\\
   U^{0} \ar[r] &[\Ve^{0}/\Hg]
   }
   \]
where the top row $\Hg \times U^{0} \arr \Ve^{0}$ is the composite of the embedding $\Hg \times U^{0} \subseteq \Hg \times \Ve^{0}$ followed by the action $\Hg\times \Ve^{0} \arr \Ve^{0}$. Since the morphism $\Ve^{0} \arr [\Ve^{0}/\Hg]$ is faithfully flat, it is enough to show that $\Hg \times U^{0} \arr \Ve^{0}$ is flat; and because $\Hg \times U^{0}$ and $\Ve^{0}$ are smooth and irreducible of dimension $16$ and $8$ respectively, it is enough to show that the codimension of the fibers is $8$ everywhere. 

We will use the formula for the action of $\Hg$ on $\Ve$ given in \ref{subsec:thickened-triple}, which is as follows:
   \begin{align*}
   (u, A + \epsilon B)\cdot \bigl(f(x) + \epsilon g(x)\bigr)
   &= \det(A + \epsilon B)^{-1}f(xA)\\
   &\qquad + \epsilon u\det(A + \epsilon B)^{-1}
   \bigl(xB \cdot \rJ_{f}(xA)  + g(xA)\bigr)\bigr).
   \end{align*}

From this formula it easy to see that $U^{0}$ does not contain any orbit. In fact, suppose that $f+\epsilon g$ is in $U^{0}$. If $g \neq 0$, there exists a matrix $A \in \GL_{2} \subseteq \Hg$ that carries $(0,1)$ into a point of $\AA^{2}$ that is not a zero of $g$, and this will carry $f+\epsilon g$ into an element of $\Ve^{0} \setminus U^{0}$. If $g = 0$, then $f \neq 0$, and since the characteristic of $k$ is different from $3$ we have $\rJ_{f}(x) \neq 0$. We can find $B \in \rM_{2}$ such that $xB \cdot \rJ_{f}(xA) \neq 0$; and the from the formula above it follows that $(\rI_{2} + \epsilon B)\cdot f = \widetilde{f} + \epsilon\widetilde{g}$ with $\widetilde{g} \neq 0$. If $\widetilde{g}(0,1) \neq 0$ we are done, otherwise we apply a matrix $A \in \GL_{2}$ as in the previous case.

The fiber of $\Hg \times U^{0} \arr \Ve^{0}$ over an element $\theta$ of $\Ve^{0}$ is the inverse image of $U^{0}$ under the morphism $\Hg \arr \Ve$ defined by $h \arrto h \cdot\theta$; since its image is not contained in $U^{0}$, by what we just said, this inverse image is a hypersurface in $\Hg$, hence it is equidimensional of dimension~$8$, as claimed.
\end{proof}

Consider the morphism
   \begin{align*}
   \beta\colon (\AA^2 \setminus \{0\}) \times \Ve &\longrightarrow \AA^3,\\
   (p, f + \epsilon g) & \longmapsto 
   \bigl(g(p), \rJ_{f}(p)\trans\bigr).
   \end{align*}

\begin{lemma}\label{lem:beta-smooth}
The morphism $\beta$ above is smooth, with irreducible geometric fibers.
\end{lemma}

\begin{proof}
We can factor $\beta$ as 
  \begin{alignat*}{3}
  (\AA^2 \setminus \{0\}) \times \Ve & {}\larr{} &
  (\AA^{2}\setminus\{0\})\times\AA^{3} & {}\larr {}\AA^{3}\\
  (p, f+\epsilon g) &\larrto & \bigl(p, g(p), \rJ_{f}(p)\trans\bigr),\\
  &&(p, v)&{}\larrto {}  v\,.
  \end{alignat*}
The first morphism is a homomorphism of vector bundles on $\AA^{2}\setminus\{0\}$, which is easily seen to be surjective; hence it is smooth with irreducible geometric fibers. The second is evidently smooth with irreducible geometric fibers. The result follows.
\end{proof}

Consider the closed subscheme $\beta^{-1}(0) \subseteq \AA^2 \setminus \{0\} \times \Ve$. This is smooth and irreducible, because of the Lemma above, and is invariant under the action of $\gm$ on $\AA^2 \setminus \{0\} \times \Ve$ by multiplication on the factor $\AA^2 \setminus \{0\}$; hence it descends to a smooth geometrically irreducible closed subscheme
   \[
   \widetilde{W} \subseteq \PP^1 \times \Ve.
   \]
We denote by $W\subseteq \Ve$ the image of $\widetilde{W}$ in $\Ve$, with its reduced scheme structure.

\begin{lemma}\label{lem:birational-1}
The projection $\widetilde{W} \arr W$ is a birational isomorphism.
\end{lemma}

\begin{proof}
The projection above is proper, and $\widetilde{W}$ is irreducible. To verify that it is birational, it is enough to check that there is a geometric point of $W$, whose inverse image in $\widetilde{W}$ is a reduced point.

Let $f \in V$ be a form of degree~$3$ with a double root $p \in \PP^{1}$ and a simple root $q \neq p$, and $g\in V$ a form with a simple root in $p$. Then $(p, f+\epsilon g)$ is the only point of $\widetilde{W}$ lying over $f + \epsilon g$, and it is easy to verify that it appears with multiplicity~$1$ in the scheme-theoretic inverse image of $f+ \epsilon g$.
\end{proof}

Notice that the group $\Hg$ acts naturally on $\PP^{1}$ through the projection $\Hg \arr \GL_{2}$ (the action of $\GL_{2}$ on $\PP^{1}$ descends from the action of $\GL_{2}$ on $\AA^{2}$ described as $(A, p) \arrto pA^{-1}$). It is easy to see that $\widetilde{W}$ is invariant under this action; hence $W \subseteq \Ve$ is an $\Hg$ invariant subscheme of $\Ve$.

Consider an object $(C \to P \to \spec k, q)$ of $\cP'_{g}(\spec k)$. This is given by a conic $P \arr \spec k$, a rational point $q \in P(k)$ and a triple cover $C \arr P$. 

\begin{lemma}\label{lem:singular-points}
The curve $C$ has a singular point lying over $q$ if and only if $(C \to P \to \spec k, q)$ is in $\Lambda^{-1}[W/\Hg]$.
\end{lemma}

\begin{proof}
We may assume that $k$ is algebraically closed. 

Consider the rank~$2$ vector bundle $E$ on $P$, with a section $\alpha \in H^0(P,\symthree{E})$, corresponding to the cover $C \arr P$. Choose a non-zero tangent vector to $P$ at $q$, corresponding to an embedding $\spec k[\epsilon] \subseteq P$, and a basis for the restriction of $E$ to $\spec k[\epsilon]$. Then the restriction of $\alpha$ to $\spec k[\epsilon]$ gives an element $f + \epsilon g$ of $\Ve$.

We need to consider two cases.

First, assume that $f = 0$. In this case the fiber of $C$ over $P$ is isomorphic to $\spec k[x,y]/\bigl((x,y)^{2}\bigr)$, so it is not Cohen--Macaulay, and $C$ must be singular. On the other hand $\epsilon g$ is in $W$, because it is the image of $(p, \epsilon g)$, where $p$ is a root of $g$.

Next suppose that $f \neq 0$. In this case the inverse image $D$ of $\spec k[\epsilon] \subseteq P$ in $C$ is the closed subscheme of $\PP^{1}_{k[\epsilon]}$ defined by $f + \epsilon g$, considered as a section of $\cO_{\PP^{1}_{k[\epsilon]}}(3)$, by Theorem~\refall{thm:describe-triple}{3}. A point of $D$ is a zero of $f$ in $\PP^{1}(k)$, and the tangent space of $C$ at $p$ coincides with the tangent space of $D$ at $p$. We have that $f + \epsilon g$ in $W$ if and only if $f$ has a zero $p$ in $\PP^{1}$ such that $\rJ_{f}(p) = g(p) = 0$. But it is immediate to check that the condition $\rJ_{f}(p) = g(p) = 0$ is equivalent to the statement that the tangent space of $D$ at $p$ is $2$-dimensional, and the Lemma follows.
\end{proof}

Now, denote by $\cS_{g}$ the complement of $\T$ in $\OT$, with its reduced structure.

\begin{lemma}\label{lem:properties-W}
Assume that $g > 1$.
\begin{enumeratea}

\item The inverse image $\Lambda^{-1}[W/ \Hg]$ is geometrically integral of codimension~$2$.

\item The image of $\Lambda^{-1}[W/ \Hg]$ under the the projection $\cP'_{g} \arr \OT$ coincides with $\cS_{g}$.

\item The morphism $\Lambda^{-1}[W/ \Hg] \arr \cS_{g}$ is a birational isomorphism.

\end{enumeratea}
\end{lemma}

\begin{proof}
Part~(a) follows from Corollary~\ref{cor:integral-inverse-image} and the fact that $W$ is integral of codimension~2. Part~(b) follows from Lemma~\ref{lem:singular-points}.

To prove part~(c) it is enough to show that there is a geometric point of $\cS_{g}$ whose inverse image in $\Lambda^{-1}[W/ \Hg]$ consists of a single reduced point. Let $C \arr \PP^{1}$ be a triple cover defined over an algebraically closed extension $\Omega$ of $k$, where $C$ is a reduced connected curve of genus~$g$ with a single node over $(1:0) \in \PP^{1}$ (it is easy to see that this exists). Then the inverse image of $\Lambda^{-1}[W/ \Hg]$ in $\PP^{1} = \spec\Omega\times_{\OT}\cP'_{g}$ consists of the single point $(1:0)$; we need to show that it appears with the reduced scheme structure. Let $t$ be a local coordinate around $(1:0)$. There exists an étale neighborhood $U \arr \PP^{1}$ of $(1:0)$, such that the pullback $C_{U}$ is isomorphic to the subscheme of $\PP^{1}_{U}$ defined by the equation $x_{1}^{2}x_{2} - t^{2}x_{2}^{3} = 0$. To calculate the composite $U \arr \PP^{1} \arr [\Ve/\Hg]$, which is associated with the embedding $C_{U} \subseteq \PP^{1}_{U}$ given above, we notice the equality
   \[
   x_{1}^{2}x_{2} - (t + \epsilon)^{2}x_{2}^{3}
      = (x_{1}^{2}x_{2} - t^{2}x_{2}^{3})
      +\epsilon 2tx_{2}^{3}
   \]
from which it follows that the morphism $U \arr [\Ve/\Hg]$ factors through the morphism $U \arr \Ve$ defined by $t \arrto (x_{1}^{2}x_{2} - t^{2}x_{2}^{3}) + \epsilon 2tx_{2}^{3}$. The subscheme $\widetilde{W}\times_{\PP^{1}\times\Ve}\PP^{1}_{U}$ is defined by setting the partial derivatives of $x_{1}^{2}x_{2} - t^{2}x_{2}^{3}$ and $2tx_{2}^{3}$ to $0$; this defines a unique reduced point of $\PP^{1}_{U}$. This implies that $W\times_{\PP^{1}\times\Ve}\PP^{1}_{U}$ also consists of a unique reduced point, and completes the proof.
\end{proof}

\section{The class of singular coverings and $\pic \T$}\label{secpic}

Now we use the description of the stack $\OT$ given in Theorem \ref{quotient} as a quotient $[X_{g}/\Gamma_{g}]$ to compute the Picard group of $\T$. We will assume $g > 1$. We will use the theory of \cite{edieit} without comments. If $G$ is a linear algebraic group over $k$ acting on a smooth $G$-scheme $X$, there are equivariant Chow rings $\rA^{*}_{G}(X)$; if $E$ is a $G$-equivariant vector bundle on $X$, then there will be Chern classes $\ch_{i}(E) \in \rA^{i}_{G}(X)$. In particular we set $\rA^{*}_{G} \eqdef \rA^{*}_{G}(\spec k)$; if $\phi\colon G \arr \GL_{n}$ is a homomorphism, we can interpret it as an equivariant vector bundle on $\spec k$, so it will have Chern classes $\ch_{i}(\phi) \in \rA^{i}_{G}$.

The Picard group $\pic\OT$ is the Picard group of equivariant invertible sheaves on $X_{g}$. By \cite[Proposition~18]{edieit}, the first Chern class homomorphism induces an isomorphism $\pic\OT \arr \rA^{1}_{\Gamma_{g}}(X_{g})$, where $\rA^{1}_{\Gamma_{g}}(X_{g})$ is the part of degree~$1$ of the equivariant Chow ring $\rA^{*}_{\Gamma_{g}}(X_{g})$. Since $X_{g}$ is a vector bundle over $\widehat{\Omega}_{g}$, we have $\rA^{1}_{\Gamma_{g}}(X_{g}) \simeq \rA^{1}_{\Gamma_{g}}(\widehat{\Omega}_{g})$, from the homotopy invariance of Chow groups. From the localization sequence
   \[
   \rA^{*}_{\Gamma_{g}}(\rM_{2, g+2} \setminus \widehat{\Omega}_{g})
   \arr \rA^{*}_{\Gamma_{g}}(\rM_{2, g+2})
   \arr \rA^{*}_{\Gamma_{g}}(\widehat{\Omega}_{g})
   \arr 0
   \]
and from Propositions \ref{prop:codimension-1} and \ref{prop:codimension-2} we see that the restriction homomorphism $\rA^{1}_{\Gamma_{g}}(\rM_{2, g+2}) \arr \rA^{1}_{\Gamma_{g}}(\widehat{\Omega}_{g})$ is an isomorphism. Furthermore, again from homotopy invariance we have that $\rA^{1}_{\Gamma_{g}}(\rM_{2, g+2})$ is isomorphic to $\rA^{1}_{\Gamma_{g}} \eqdef \rA^{1}_{\Gamma_{g}}(\spec k)$. Applying once again \cite[Proposition~18]{edieit}, we see that $\rA^{1}_{\Gamma_{g}}$ is isomorphic to the group of characters $\widehat{\Gamma}_{g}$ of $\Gamma_{g}$, which is easily seen to be isomorphic to $\ZZ \oplus \ZZ$. 

On the other hand we are interested in the Picard group of $\T$, which equals $\OT \setminus \cS_{g}$ for a certain closed integral substack $\cS_{g} \subseteq \OT$ of codimension~$1$, by Lemma~\ref{lem:properties-W}. The substack $\cS_{g} \subseteq \OT \simeq [X_{g}/\Hg]$ is of the form $[Y_{g}/\Hg]$, where $Y_{g}$ is an integral $\Hg$-invariant hypersurface in $X_{g}$. From the localization sequence
   \[
   \xymatrix@R6pt{
   \rA^{0}_{\Gamma_{g}}(Y_{g}) \ar[r]\ar@{ = }[d]&
   \rA^{1}_{\Gamma_{g}}(X_{g}) \ar[r]\ar@{ = }[d]&
   \rA^{1}_{\Gamma_{g}}(X_{g}\setminus Y_{g}) \ar[r]\ar@{ = }[d]&
   0\\
   \ZZ & \rA^{1}(\OT)& \rA^{1}(\T)
   }
   \]
we deduce that $\rA^{1}(\T)$ is the quotient of $\rA^{1}_{\Gamma_{g}} \simeq \ZZ \oplus \ZZ$ by the subgroup generated by the class of $Y_{g}$ in $\rA^{1}_{\Gamma} = \rA^{1}_{\Gamma_{g}}(X_{g})$. We need to compute this class.

It is convenient to work with the group
   \[
   \Delta_{g} \eqdef \GL_{g+4}\times \GL_{g+2} \times \GL_{2}.
   \]
Recall that $\Gamma_{g}$ is the quotient of $\Delta_{g}$ by the subgroup $\gm \subseteq \Delta_{g}$, where $\gm$ is embedded into $\Delta_{g}$ by the homomorphism $t \arrto (\rI_{g+4}, t\rI_{g+2}, t^{-1}\rI_{2})$. The character group $\widehat{\Gamma}_{g}$ is the kernel of the restriction homomorphism $\widehat{\Delta}_{g} \arr \widehat{\GG}_{\mathrm{m}}$. Let us translate this in terms of Chow groups.

The Chow ring $\rA^{*}_{\Delta_{g}}$ is a polynomial ring
   \[
   \ZZ[\delta_{1}, \dots, \delta_{g+4}, \gamma_{1}, \dots, \gamma_{g+2}, 
   \sigma_{1}, \sigma_{2}],
   \]
where the $\delta_{i}$ are the Chern classes of the representation given by the first projection $\Delta_{g} \arr \GL_{g+4}$, the $\gamma_{i}$ the Chern classes of the second projection $\Delta_{g} \arr \GL_{g+2}$, and the $\sigma_{i}$ of the third projection $\Delta_{g} \arr \GL_{2}$. Then $\rA^{1}_{\Delta_{g}}$ is the free abelian group generated by $\delta_{1}$, $\gamma_{1}$ and $\sigma_{1}$. Then the natural homomorphism $\rA^{1}_{\Gamma_{g}} \arr \rA^{1}_{\Delta_{g}}$ identifies $\rA^{1}_{\Gamma_{g}}$ with the kernel of the homomorphism $\rA^{1}_{\Delta_{g}} \arr \rA^{1}_{\gm}$. If we denote by $\tau$ the canonical generator of $\rA^{1}_{\gm}$, the homomorphism $\rA^{1}_{\Delta_{g}} \arr \rA^{1}_{\gm}$ sends $\delta_{1}$ to $0$, $\gamma_{1}$ to $(g+2)\tau$ and $\sigma_{1}$ into $-2\tau$; hence the kernel is generated by $\delta_{1}$ and $q_{1}$, where 
   \begin{equation}\label{eq:q1}
   q_{1} \eqdef
   \begin{cases}
   (g+2)\sigma_{1} + 2\gamma_{1} & \text{if $g$ is odd, and}\\
   \frac{\displaystyle g+2}{\displaystyle 2}\sigma_{1} + \gamma_{1} 
   & \text{if $g$ is even.}
   \end{cases}
   \end{equation}
Thus we need to compute the class of $Y_{g}$ in $\rA^{1}_{\Delta_{g}}$ and express it in terms of the generators $\delta_{1}$, $q_{1}$.

\subsection*{The class of $W$ in $\rA^{2}_{\Hg}$.}

As a first step, we compute the class of $W$ in $\rA^{2}_{\Hg}(\Ve) = \rA^{2}_{\Hg}$.
According to \cite[Lemma~2.1]{vismol} there is a canonical isomorphism of graded rings $\rA^*_{\Hg}(\Ve)\cong \rA^*_{\gm\times\GL_{2}}(\Ve)$ induced by the embedding $\gm \times \GL_{2} \subseteq \Hg$;
the action of $\gm \times \GL_{2}$ on $\Ve$ is given by the formula
   \[
   (\alpha, A)(f+\epsilon g)(x)
   = \det(A)^{-1}\bigl(f(xA) + \epsilon \alpha g(xA)\bigr).
   \]
The ring $\rA^*_{\gm\times\GL_{2}}$ is a polynomial ring $\ZZ[\nu_{1}, c_{1}, c_{2}]$, where $\nu_{1}$ is the first Chern class of the representation given by the projection $\gm \times \GL_{2}\arr \gm$, while $c_{1}$ and $c_{2}$ are the Chern classes of the second projection $\gm \times \GL_{2} \arr \GL_{2}$.

By Lemma~\ref{lem:birational-1}, the class of $W$ in $\rA^*_{\gm \times \GL_{2}}(\Ve)$ is the pushforward of the class of $\widetilde{W}$ in $\rA^*_{\gm \times \GL_{2}}(\PP^{1}\times\Ve)$. We can also write $\rA^*_{\gm \times \GL_{2}}(\PP^{1}\times\Ve)$ as
   \[
   \rA^*_{\gm \times\gm \times \GL_{2}}
   \bigl((\AA^{2}\setminus\{0\})\times\Ve\bigr),
   \]
where the action of $\gm \times\gm \times \GL_{2}$ is defined as
   \[
   (\lambda, \alpha, A)\bigl(p, (f+\epsilon g)(x)\bigr)
   = \bigl(\lambda pA^{-1},
      \det(A)^{-1}(f(xA) + \epsilon \alpha g(xA))\bigr).
   \]
(recall that we are writing the vectors in $\AA^{2}$ as row vectors). 

Consider the morphism
   \begin{align*}
   \beta\colon (\AA^2 \setminus \{0\}) \times \Ve &\longrightarrow \AA^3\\
   (p, f + \epsilon g) & \longmapsto 
   \bigl(g(p), \rJ_{f}(p)\trans\bigr)
   \end{align*}
of the preceding section. We have
   \begin{alignat*}{3}
   \beta\bigl((\lambda,\alpha,A)\cdot\bigl(p,f(x)+\epsilon g(x)\bigr)\bigr) &=
   \beta\bigl(\lambda pA^{-1},\det(A)^{-1}(f(xA)
      + \epsilon \alpha g(xA))\bigr)\\
   &=\bigl(\det(A)^{-1}\alpha g(\lambda p),
      \rJ_{\det(A)^{-1}f(xA)}(\lambda pA^{-1})\trans\bigr)\\
   &=\bigl(\det(A)^{-1}\alpha\lambda^{3} g(p),
      \det(A)^{-1}\lambda^{2}\rJ_{f(xA)}(pA^{-1})\trans\bigr)\\
   &=\bigl(\det(A)^{-1}\alpha\lambda^{3} g(p),
      \det(A)^{-1}\lambda^{2} (A\rJ_{f}(p))\trans\bigr)\\
   &=\bigl(\det(A)^{-1}\alpha\lambda^{3} g(p),
      \det(A)^{-1}\lambda^{2} \rJ_{f}(p)\trans A\trans\bigr) 
   \end{alignat*}
where the penultimate equality is obtained by the chain rule. Thus we define a linear  action of $\gm \times \gm \times \GL_{2}$ on $\AA^{3} = \AA^{1} \oplus \AA^{2}$ by the rule
   \[
   (\lambda, \alpha, A)(v, w) =
   \bigl(\det(A)^{-1}\alpha\lambda^{3} v,
      \det(A)^{-1}\lambda^{2} wA\trans\bigr);
   \]
this makes the morphism $\beta$ into a $\gm \times \gm \times \GL_{2}$-equivariant morphism.

Denote by $\nu_{1}$, $c_{1}$ and $c_{2}$ the pullbacks to $\rA^*_{\gm \times \gm \times\GL_{2}}$ of the corresponding classes in $\rA^*_{\gm \times \GL_{2}}$, and by $\mu_{1}$ the first Chern class of the representation $\gm \times \gm \times \GL_{2} \arr \gm$ given by the first projection. The class of $\{0\}$ in $\rA^{*}_{\gm \times \gm \times \GL_{2}}(\AA^{3}) = \rA^{*}_{\gm \times \gm \times \GL_{2}}$ is the third Chern class of $\AA^{3} = \AA^{1} \oplus \AA^{2}$. If $a_{1}$ $a_{2}$ are the Chern roots of the projection $\gm \times \gm \times \GL_{2} \arr \GL_{2}$, the second Chern class of $\AA^{2}$ is
   \[
   (a_{1}-c_{1}+2\mu_{1})(a_{2}-c_{1}+2\mu_{1}) =
   c_{2} - 2c_{1}\mu_{1} + 4\mu_{1}^{2}.
   \]
The action of $(\lambda, \alpha, A) \in \gm \times \gm \times \GL_{2}$ on the first factor $\AA^{1}$ is multiplication by $\det(A)^{-1}\alpha\lambda^{3}$; hence the third Chern class of $\AA^{3}$ is
\begin{equation}\label{eq:c3A3}
(c_{2} - 2c_{1}\mu_{1} + 4\mu_{1}^{2})(-c_{1}+\nu_{1}+3\mu_{1}).
\end{equation}   
The ring $\rA^{*}_{\gm \times \gm \times \GL_{2}}\bigl((\AA^{2} \setminus \{0\})\times \Ve\bigr)$ is the quotient of $\rA^*_{\gm \times \gm \times \GL_{2}}$ by the ideal generated by the second Chern class of $\AA^{2}$, considered as a linear representation of $\gm\times \GL_{2}$ with the action $(\lambda, A) \cdot p = \lambda pA^{-1}$, by \cite[Lemma~2.2]{vismol}. The second Chern class of $\AA^{2}$ is 
   \[
   (-a_{1} + \mu_{1})(-a_{2} + \mu_{1}) = \mu_{1}^{2} - c_{1}\mu_{1} + c_{2};
   \]
hence we have
   \[
   \rA^*_{\gm \times \gm \times \GL_{2}}
      \bigl((\AA^{2} \setminus \{0\})\times \Ve\bigr))
   = \ZZ[\mu_{1}, \nu_{1}, c_{1}, c_{2}]/(\mu_{1}^{2} - c_{1}\mu_{1} + c_{2}).
   \]
The class of the pullback $\beta^{-1}(0)$ in $\rA^*_{\gm \times \gm \times \GL_{2}}\bigl((\AA^{2} \setminus \{0\})\times \Ve\bigr))$ is obtained by making the substitution $\mu_{1}^{2} \arrto c_{1}\mu_{1} - c_{2}$ in (\ref{eq:c3A3}), and equals
   \begin{equation}\label{eq:class-Wtilde}
   -3 c_1 c_2 - 3c_2\nu_1 + (4c_1^2 - 9c_2 + 2c_{1}\nu_1)\mu_1.
   \end{equation}

In the isomorphism $\rA^*_{\gm \times \gm \times \GL_{2}}\bigl((\AA^{2} \setminus \{0\})\times \Ve\bigr) \simeq \rA^*_{\gm \times \GL_{2}}(\PP^{1}\times \Ve)$ the class of $\beta^{-1}(0)$ corresponds to the class of $\widetilde{W}$, and $\mu_{1}$ corresponds to the first Chern class of $\cO_{\PP^{1}}(1)$. By Lemma~\ref{lem:birational-1}, the class of $W$ in $\rA^*_{\gm\times\GL_{2}}(\Ve)$ is the pushforward of the class of $\widetilde{W}$. Since the pushforward of $1$ is $0$ and the pushforward of $\mu_{1}$ is $1$, from the projection formula we obtain the following.

\begin{lemma}\label{lem:class-W}
The class of $W$ in $\rA^{2}_{\Hg}(\Ve) = \rA^{2}_{\Hg}$ is
   \[
   4c_1^2 - 9c_2 + 2c_{1}\nu_1\,.
   \]
\end{lemma}

\subsection*{The class of $Y_{g}$ in $\rA^{1}_{\Delta_{g}}$.}

Consider the morphism $\Lambda\colon \cP'_{g} \arr [\Ve/\Hg]$. The stack $\cP'_{g}$ is the quotient $[(\PP^{1}\times X_{g})/\Gamma_{g}]$, where $\Gamma_{g}$ acts on $X_{g}$ in the way detailed above, while the action on $\PP^{1}$ is described by the obvious homomorphism $\Gamma_{g} \arr \PGL_{2}$. Since $\cP'_{g} = [(\PP^{1}\times X_{g})/\Gamma_{g}]$ and $\OT = [X_{g}/\Gamma_{g}]$, we have a commutative diagram
   \[
   \xymatrix{
   \PP^{1}\times X_{g} \ar[r] \ar[d]&
   [(\PP^{1}\times X_{g})/\Delta_{g}] \ar[r] \ar[d]&
   \cP'_{g}\ar[r]^-{\Lambda} \ar[d]&
   [\Ve/\Hg]\\
   X_{g} \ar[r]&
   [X_{g}/\Delta_{g}] \ar[r]&
    \OT\,.
   }
   \]
Denote by $\widetilde{Y}_{g}$ the inverse image of $[W/\Hg]$ in $\PP^{1}\times X_{g}$ and by $Y_{g}$ its image in $X_{g}$. We need to compute the class of $Y_{g}$ in the equivariant Chow ring $\rA^{*}_{\Delta_{g}}(X_{g})$; Lemma~\ref{lem:properties-W} implies that this is the pushforward of the class of $\widetilde{Y}_{g}$ in $\rA^{2}_{\Delta_{g}}(\PP^{1}\times X_{g})$. 

Denote by $\xi$ the first Chern class of the sheaf $\cO(1)$ on $\PP^{1}$, considered as a $\Delta_{g}$-space with the action coming from the projection $\Delta_{g} \arr \GL_{2}$. Since $X_{g}$ is a vector bundle over an open subscheme of a representation of $\Delta_{g}$, $\rA^{*}_{\Delta_{g}}(\PP^{1}\times X_{g})$ is a quotient of $\rA^*_{\Delta_{g}}(\PP^{1})$. Using the notation for the classes in $\rA^*_{\Delta_{g}}$ introduced at the beginning of the section, the ring $\rA^*_{\Delta_{g}}(\PP^{1})$ has the form
   \[
   \ZZ[\delta_{1}, \dots, \delta_{g+4}, \gamma_{1}, \dots, \gamma_{g+2},
   \sigma_{1}, \sigma_{2}, \xi]/(\xi^{2} + \sigma_{1}\xi + c_{2}).
   \]

The class of $\widetilde{Y}_{g}$ in $\rA^{*}_{\Delta_{g}}(\PP^{1}\times X_{g})$ is the pullback of the class of $W$ in $\rA^{*}_{\Hg}(\Ve)$. Now, by Proposition~\ref{isoquot} and the construction of the morphism $\Lambda$, the pullback of the canonical invertible sheaf on $[\Ve/\Hg]$ is $\Omega^{1}_{\cP'_{g}/\OT}$, while the pullback of the locally free of rank~$2$ that comes from the projection $\Hg \arr \GL_{2}$ is the tautological sheaf $\cE'_{g}$. This gives a formula for the pullback $\rA^{*}_{\Hg} = \rA^{*}_{\Hg}(\Ve) \arr \rA^{*}(\cP'_{g})$: the  class $\nu_{1}$ goes to the first Chern class of $\Omega^{1}_{\cP'_{g}/\OT}$, while the $c_{i}$ go to the Chern classes of $\cE'_{g}$.

We need to compute the first Chern class of $\Omega^{1}_{(\PP^{1}\times X_{g})/X_{g}}$ in $\rA^{1}_{\Delta_{g}}(\PP^{1}\times X_{g})$, which is the pullback of the first Chern class of $\Omega^{1}_{\PP^{1}/\spec k}$ in $\rA^{1}_{\Delta_{g}}(\PP^{1})$. Since $\PP^{1}$ is the projective space of lines in $\AA^{2}$ (its projectivization in the classical sense, dual to Grothendieck's), and the Chern classes of $\AA^{2}$ considered as a equivariant vector bundle on $\spec k$ are $\sigma_{1}$ and $\sigma_{2}$, the class $\xi \eqdef \rmc_{1}(\AA^{2})$ satisfies the relation
\begin{equation}\label{eq:class-xi}
\xi^{2} + \sigma_{1}\xi + \sigma_{2} = 0\in \rA^{2}(\PP^{1})\,.
\end{equation}
Furthermore, it follows from the Euler sequence that
\begin{equation}\label{eq:class-omega}
\rmc_{1}(\Omega^{1}_{\PP^{1}/\spec k}) = -2\xi - \sigma_{1} \in \rA^{2}(\PP^{1})\,.
\end{equation}

To compute the Chern classes of $\cE'_{g}$, recall from Section~\ref{sec:glgen-conic} that by construction this fits into a $\Delta_{g}$-equivariant exact sequence
   \[
   0 \arr \cO_{\PP^1_{X_{g}}}(-1)^{g+2} \arr \cO_{\PP^1_{X_{g}}}^{g+4} \arr \cE_g' \arr 0\,.
   \]
In the calculations that follow we only need to record the elements of degree at most~$2$. We will write $\equiv$ to mean ``equal up to terms of degree $> 2$''.

The Chern classes of $\cO_{\PP^1_{X_{g}}}^{g+2}$ are the $\gamma_{i}$; denote by $a_{1}$, \dots, $a_{g+2}$ its Chern roots.  We have
   \begin{align*}
   \ch\bigl(\cO_{\PP^1_{X_{g}}}(-1)^{g+2}\bigr) &=
   (1 + a_{1} -\xi) \dots (1 + a_{g+2} -\xi)\\
   &\equiv 1 + \gamma_{1} - (g+2)\xi 
      +\gamma_{2} - (g+1)\gamma_{1}\xi + \frac{(g+1)(g+2)}{2}\xi^{2}
   \end{align*}
and
   \begin{align*}
   \ch\bigl(\cO_{\PP^1_{X_{g}}}(-1)^{g+2}\bigr)^{-1} &\equiv
   1 - \bigl(\gamma_{1} - (g+2)\xi\bigr)
      + \bigl(\gamma_{1} - (g+2)\xi\bigr)^{2}\\
      &\qquad - \Bigl(\gamma_{2} - (g+1)\gamma_{1}\xi +
         \frac{(g+1)(g+2)}{2}\xi^{2}\Bigr)\\
   &=
   1 - \gamma_{1} + (g+2)\xi
   + \gamma_{1}^{2} - \gamma_{2}
   -(g+3)\gamma_{1}\xi
   +\frac{(g+2)(g+3)}{2}\xi^{2}\,;
   \end{align*}
from this we obtain that
   \begin{align*}
   \ch(\cE'_{g}) &= \ch\bigl(\cO_{\PP^1_{X_{g}}}(-1)^{g+2}\bigr)^{-1}
      \ch\bigl(\cO_{\PP^1_{X_{g}}}^{g+4}\bigr)\\
   &\equiv \bigl(1 - \gamma_{1} + (g+2)\xi
      + \gamma_{1}^{2} - \gamma_{2}
      -(g+3)\gamma_{1}\xi
      +\frac{(g+2)(g+3)}{2}\xi^{2}\bigr)\\
   &\qquad(1 + \delta_{1} + \delta_{2})\\
   &\equiv 1 + \bigl(\delta_{1} - \gamma_{1} + (g+2)\xi\bigr)\\
   &\qquad+ \delta_{2} - \delta_{1}\gamma_{1}
      + \gamma_{1}^{2} - \gamma_{2}
      +(g+2)\delta_{1}\xi
      -(g+3)\gamma_{1}\xi
      +\frac{(g+2)(g+3)}{2}\xi^{2}.
   \end{align*}

Using the relation $\xi^{2} + \sigma_{1}\xi + \sigma_{2} = 0$, we can rewrite this as
   \begin{align*}
   \ch(\cE'_{g}) &= 
   1 + \delta_{1} - \gamma_{1} + (g+2)\xi\\
   &\qquad+ \delta_{2} - \gamma_{1}\delta_{1}+ \gamma_{1}^{2} - \gamma_{2}
   -\frac{(g+2)(g+3)}{2}\sigma_{2}\\
   &\qquad + \Bigl(
   (g+2)\delta_{1} - (g+3)\gamma_{1} - \frac{(g+2)(g+3)}{2}\sigma_{1}
   \Bigr)\xi.
   \end{align*}

So, the pullback $\rA_{\Hg}[\Ve/\Hg] \arr \rA_{\Delta_{g}}(\PP^{1}\times X_{g})$ acts as follows:
   \begin{align*}
   \nu_{1} &\larrto -2\xi - \sigma_{1}\\
   c_{1}&\larrto \delta_{1} - \gamma_{1} + (g+2)\xi\\
   c_{2}&\larrto 
   \delta_{2} - \gamma_{1}\delta_{1}+ \gamma_{1}^{2} - \gamma_{2}
   -\frac{(g+2)(g+3)}{2}\sigma_{2}\\
   &\qquad + \bigl(
   (g+2)\delta_{1} - (g+3)\gamma_{1} - \frac{(g+2)(g+3)}{2}\sigma_{1}
   \bigr)\xi\,;
   \end{align*}
hence the class of $\widetilde{Y}_{g}$ in $\rA^{2}_{\Delta_{g}}(\PP^{1}\times X_{g})$, which is the pullback of the class of $W$ in $\rA^2_{\Hg}(\Ve)$, is obtained by applying the substitution above to the equation for the class of $W$ given in Lemma~\ref{lem:class-W}, and then getting rid of the term $\xi^{2}$ via the substitution $\xi^{2} \arrto -\sigma_{1}\xi - \sigma_{2}$. After some calculations, one obtains that the class of $\widetilde{Y}_{g}$ is
   \begin{align*}
   &-5\gamma_{1}^{2} + \gamma_{1}\delta_{1} + 2\gamma_{1}\sigma_{1}
      +9\gamma_{2} + 4\delta_{1}^{2} + 2\delta_{1}\sigma_{1} -9\delta_{2}
      +\frac{(g+2)(19g+73)}{2}\sigma_{2}\\
   &\quad +\Bigl(
   (g+15)\gamma_{1} - (g + 6)\delta_{1} -\frac{(g+2)(g + 15)}{2}\sigma_{1}
   \Bigr)\xi.
   \end{align*}
Finally we use Lemma~\ref{lem:properties-W}~(c) to notice that the class of $Y_{g}$ in $\rA^{1}_{\Delta_{g}}(X_{g})$ is the pushforward of the class of $\widetilde{Y_{g}}$. By the projection formula we obtain that this is the coefficient of $\xi$ in the expression above for the class of $\widetilde{Y_{g}}$, which is
   \[
   [Y_{g}] = 
   (g+15)\gamma_{1} - (g + 6)\delta_{1} -\frac{(g+2)(g + 15)}{2}\sigma_{1}\,.
   \]
This must lie in the subgroup $\rA^{1}_{\Gamma_{g}}(X_{g}) \subseteq  \rA^{1}_{\Delta_{g}}(X_{g})$, which, as we know, is generated by $\delta_{1}$ and $q_{1}$ (the expression for $q_{1}$ is given in (\ref{eq:q1})). With considerable relief we notice that it is in fact so, and
   \[
   [Y_{g}] = 
   \begin{cases}
   -(g + 6)\delta_{1}
      + \frac{\displaystyle g+15}{\displaystyle\mathstrut 2}q_{1}
   & \text{if $g$ is odd, and}\\
   -(g + 6)\delta_{1} - (g + 15)q_{1}
   &\text{if $g$ is even.}
   \end{cases}
   \]
So
   \[
   \pic\T \simeq
   \begin{cases}
   (\ZZ\oplus\ZZ)/\generate{\bigl(-(g + 6), 
      \frac{\displaystyle (g + 15)}{\displaystyle\mathstrut 2}\bigr)}
   & \text{if $g$ is odd, and}\\
   (\ZZ\oplus\ZZ)/\generate{\bigl(-(g + 6), g + 15\bigr)}
   &\text{if $g$ is even.}
   \end{cases}
   \]
An elementary calculation ends the proof of Theorem~\ref{thm:picard}.


\bibliographystyle{amsalpha}
\bibliography{bibtrig}

\end{document}